\documentclass[11pt,letterpaper]{amsart}
\usepackage{amssymb, amscd, amsfonts, euler, epic, xy}
\xyoption{all}
\usepackage{charter}

\def\CC{{\mathbb C}}

\def\FF{{\mathbb F}}

\def\PP{{\mathbb P}}
\def\QQ{{\mathbb Q}}
\def\RR{{\mathbb R}}

\def\ZZ{{\mathbb Z}}

\def\hfrak{{\mathfrak h}}
\def\mfrak{{\mathfrak m}}
\def\tfrak{{\mathfrak t}}

\def\half{\tfrac{1}{2}}
\def\G{{\Gamma}}

\def\Det{\mathit{Det}}

\def\spec{{\rm Spec}}

\def\reg{{\rm reg}}

\def\bs{\backslash}

\def\pt{{\bullet}}
\def\eps{\epsilon}

\def\coble{{\mathcal C}}

\def\Ccal{{\mathcal C}}

\def\Ical{{\mathcal I}}

\def\Lcal{{\mathcal L}}

\def\Mcal{{\mathcal M}}

\def\Ocal{{\mathcal O}}

\def\Scal{{\mathcal S}}

\def\la{\langle}
\def\ra{\rangle}

\newcommand\bl{\operatorname{Bl}}

\newcommand\Gr{\operatorname{Gr}}
\newcommand\Hom{\operatorname{Hom}}

\newcommand\pic{\operatorname{Pic}}

\newcommand\proj{\operatorname{Proj}}
\newcommand\res{\operatorname{Res}}
\newcommand\sym{\operatorname{Sym}}

\newcommand\GL{\operatorname{GL}}

\newcommand\PGL{\operatorname{PGL}}

\newcommand\SL{\operatorname{SL}}

\newcommand\Sp{\operatorname{Sp}}

\newcommand\sign{\operatorname{sign}}

\newtheorem{theorem}{Theorem}[section]
\newtheorem{lemma}[theorem]{Lemma}
\newtheorem{proposition}[theorem]{Proposition}

\newtheorem{lemmadef}[theorem]{Lemma-Definition}

\newtheorem{corollary}[theorem]{Corollary}
\newtheorem{conjecture}[theorem]{Conjecture}

\theoremstyle{definition}
\newtheorem{definition}[theorem]{Definition}

\theoremstyle{remark}
\newtheorem{remark}[theorem]{Remark}

\newtheorem{question}[theorem]{Question}

\def\DynkinEEEE#1#2#3#4#5#6#7#8
   {\begin{picture}(124, 55)%
   \put(5, 40){\circle{4}}%
   \put(7, 40){\line(1, 0){28}}%
   \put(37, 40){\circle{4}}%
   \put(39, 40){\line(1, 0){28}}%
   \put(69, 40){\circle{4}}%
   \put(71, 40){\line(1, 0){28}}%
   \put(101, 40){\circle{4}}%
    \dashline{2}(103, 40)(141, 40)%
    \put(143, 40){\circle{4}}%
   \put(145, 40){\line(1, 0){28}}%
   \put(175, 40){\circle{4}}%
   \put(69, 38){\line(0, -11){28}}%
   \put(69, 8){\circle{4}}%
   \put(5, 46){\makebox(0, 0)[b]{$#2$}}%
   \put(37, 46){\makebox(0, 0)[b]{$#3$}}%
   \put(69, 46){\makebox(0, 0)[b]{$#4$}}%
   \put(101, 46){\makebox(0, 0)[b]{$#5$}}%
   \put(143, 46){\makebox(0, 0)[b]{$#6$}}%
   \put(69, 5){\makebox(0, 0)[t]{$#1$}}%
    \put(175, 46){\makebox(0, 0)[lb]{$#7$}}%
\end{picture}}

\title{Del Pezzo moduli via root systems}

\author{Elisabetta Colombo}
\address{Dipartimento di Matematica, Universit\`a di Milano,
  via Saldini 50, I-20133 Milano, Italia}
 \email{elisabetta.colombo@mat.unimi.it, geemen@mat.unimi.it}
\author{Bert van Geemen}
\author{Eduard Looijenga}
\address{Mathematisch Instituut\\
Universiteit Utrecht\\
P.O.~Box 80.010, NL-3508 TA Utrecht\\
Nederland}
\email{looijeng@math.uu.nl}
\begin{document}
\begin{abstract}
Coble defined in his 1929 treatise invariants for cubic surfaces
and quartic curves. We reinterpret these in terms of the root systems
of type $E_6$ and $E_7$ that are naturally associated to these varieties, thereby
giving some of his results a more intrinsic treatment. Our discussion is uniform for
all Del Pezzo surfaces of degree 2,3,4 and 5.
\end{abstract}

\maketitle
\hfill{\emph{To Professor Yuri Manin, for his 70th birthday.}}

\section*{Introduction}
A \emph{Del Pezzo surface of degree $d$} is a smooth projective surface with semi-ample
anticanonical bundle whose class has self-intersection $d$. The degree is always
between 1 and 9 and the surface is either a quadric ($d=8$ in that case) or is obtained from
blowing up  $9-d$ points in the projective plane that satisfy a mild genericity
condition. So moduli only occur for $1\le d\le 4$. The anticanonical system
is $d$-dimensional and when $d\not=1$, it is also base point free. For $d=4$, the resulting  morphism is  
birational onto a complete intersection of two quadrics in $\PP^4$, for
$d=3$ it is birational onto a cubic surface in $\PP^3$ and for $d=2$ we get a degree
two map onto $\PP^2$ whose discriminant curve is a quartic (we will ignore the
case $d=1$ here). This image surface
(resp.\ discriminant) is smooth in case the anticanonical bundle is ample;
we then call the Del Pezzo surface a \emph{Fano surface}. Otherwise it might have simple singularities in the 
sense of Arnol'd (that have a
root system label $A$, $D$ or $E$). Conversely, every complete
intersection of two quadrics in $\PP^4$, cubic surface in $\PP^2$ or quartic curve
in $\PP^2$ with such
singularities thus arises.

We mentioned that a degree $d$ Del Pezzo surface, $d\not=8$,  is obtained from blowing up
$9-d$ points in $\PP^2$ in general position. A more precise statement is that
if we are given as many disjoint exceptional curves $E_1,\dots ,E_{9-d}$ on the Del Pezzo
surface $S$, then these can be simultaneously contracted to produce a projective plane.
So the images of these curves yield $9-d$ numbered points $p_1,\dots ,p_{9-d}$ in $\PP^2$
given up to projective equivalence. Hence every polynomial expression in terms of the
projective coordinates of these points that is invariant under $\SL(3,\CC)$ is
a `covariant' for the tuple $(S;E_1,\dots ,E_{9-d})$. Coble exhibited such
covariants for the important cases $d=2$ and $d=3$. These are in general not
covariants of $S$ itself, since the surface may have many exceptional systems
$(E_1,\dots ,E_{9-d})$. Indeed, if we assume that $S$ is Fano, then, as Manin observed, these
systems are simply transitively permuted  by a Weyl group $W$, which acts here as
a group of Cremona transformations. Therefore, this group
will act on the space of such covariants. Coble's covariants span a $W$-invariant
subspace and Coble was able to identify the $W$-action as a Cremona group
(although the Weyl group interpretation was not available to him). For $d=3$ he found
an irreducible representation of degree $10$ of a $E_6$-Weyl group and for
$d=2$ he obtained an irreducible representation of degree $15$ of a $E_7$-Weyl group.

The present paper purports to couch Coble's results in terms of a moduli space
of tuples $(S;E_1,\dots ,E_{9-d})$ as above, for which
$S$ is semistable in the sense of Geometric Invariant Theory.
This moduli space comes with an action of the Weyl group $W$. It also
carries a natural line bundle, called  the \emph{determinant bundle},
to which the $W$-action lifts: this line bundle assigns to a Del Pezzo surface
the line that is the dual of the top exterior power of the space of sections
of its anticanonical bundle. In turns out that this bundle is proportional
to the one that we use to do geometric invariant theory with (and from which our notion
of semistability originates).
We show that the Coble covariants can be quite naturally understood as sections
of this bundle  and we reprove the fact known to Coble that these sections span an irreducible 
representation of $W$.
We also show that these sections separate the points of the above moduli space
so that one might say that Coble's covariants of a stable tuple $(S;E_1,\dots ,E_{9-d})$
make up a complete set of invariants. This approach not only
covers the cases Coble considered (degree $2$ and $3$), but also the degree $4$
case and, somewhat amusingly, even the degree $5$ case, for which there are no
moduli at all. For the case of degree $3$  we also make the connection with earlier
work of Naruki and Yoshida. This allows us to conclude that the Coble covariants define a 
complete linear system
and define a closed immersion of the GIT-compactification of the moduli space of
marked cubic surfaces in a $9$-dimensional projective space.
Our results  are less complete when the degree is 2; for instance, we did not manage
to establish that the Coble covariants define a complete linear system.

We end up with a description of the GIT moduli space that is
entirely in terms of the corresponding root system. Our results lead to us to some
remarkable integrability properties of the module of $W$-invariant vector fields on the vector
space that underlies the defining (reflection) representation of $W$ and we raise the
question of whether this is a special case of a general phenomenon.

Since the appearance of Coble's book a great deal of work
on Del Pezzo moduli has seen the day.
As its sheer volume makes it impossible to give
our predecessors their fair due,
any singling out of contributions will be biased.
While keeping that in mind we nevertheless
wish to mention the influential book by Manin \cite{manin}, the Ast\'erisque volume
by Dolgachev-Ortland \cite{dolgort},
Naruki's construction of a smooth compactification
of the moduli space of marked cubic surfaces \cite{naruki},
the determination of its Chow groups in \cite{elbert},
the Lecture Note by Hunt \cite{hunt}
and Yoshida's revisit of the Coble covariants \cite{yoshida}.
The ball quotient description of the moduli space of cubic surfaces by Allcock, Carlson and Toledo \cite{act},
combined with Borcherds' theory of modular forms, led Allcock and Freitag
\cite{af} to construct an embedding of the moduli space of marked cubic surfaces, which coincides with the map given by the Coble invariants
\cite{freitag2}, \cite{geemen}.
\\
We now briefly review the organization of the paper. The first section introduces a moduli space
for marked Fano surfaces of degree $d\ge 2$ as well as the line bundle over that space that is central
to this paper, the determinant line bundle. This assigns to a Fano surface the determinant line
of the dual of the space of sections of its anticanonical sheaf (this is also the determinant of the
cohomology of its structure sheaf). We show that this line bundle can be used to obtain in a
uniform manner a compactification (by means of GIT) so the determinant bundle extends over this
compactification as an ample bundle. In Section 2 we introduce the Coble covariants and show that
they can be identified with sections of the determinant bundle. The next section expresses
these covariants purely in terms of the associated root system. In Section 4 we identify
(and discuss) the Weyl group representation spanned by the Coble covariants.  The final section investigates the separating properties of the Coble covariants, where the emphasis is on the
degree $3$ case.
\\

As we indicated, Manin's work on Del Pezzo surfaces has steered this beautiful
subject in a new direction. Although this represents only a small part of his many influential
contributions to mathematics, we find it therefore quite appropriate to dedicate
this paper to him on the occasion of his 70th birthday.

\tableofcontents

\section{Moduli spaces for marked Del Pezzo surfaces}\label{sect:gitmoduli}
We call a smooth complete surface $S$ a \emph{Del  Pezzo surface of degree $d$} if
its anticanonical bundle $\omega_S^{-1}$ is semi-ample  and $\omega_S\cdot \omega_S=d$. It is known that  then $1\le d\le 9$ and that $S$ is  isomorphic either to  a smooth quadric or to a surface obtained from successively blowing up $9-d$ points of $\PP^2$. In order that  a successive blowing up  of $9-d$ points of $\PP^2$ yields a Del  Pezzo surface  it is necessary and sufficient that we blow up on (i.e., over the strict transform of) a smooth cubic curve (which is an anticanonical divisor of $\PP^2$). This is equivalent to the apparently weaker condition that we blow up at most $3$ times on a line and at most $6$ times on a conic. It is also equivalent to the apparently stronger condition that the anticanonical
system on this surface is nonempty and has dimension $d$.

Let $S$ be a Del  Pezzo surface of degree $d$. The vector space
$V(S):=H^0(\omega_S^{-1})^*$ (which we will usually abbreviate by $V$) has dimension $d+1$. The anticanonical system defines the  (rational) \emph{anticanonical map} $S\dashrightarrow \PP(V)$. When $d\ge 2$, it  has  no base points,
so that the anticanonical map is a morphism.  For $d=1$, it has a single base point; if we blow up this point, then the anticanonical map lifts to a morphism $\tilde S\to \PP(V)$ which makes $\tilde S$ a rational elliptic surface (with  a section defined by the exceptional curve of the blowup).

Let $S\to \bar S$ contract the $(-2)$-curves on $S$
(we recall that a curve on a smooth surface is called a \emph{$(-2)$-curve} if it is a
smooth rational curve with self-intersection $-2$.) Then
$\bar S$ has rational double point singularities only and its dualizing sheaf  $\omega_{\bar S}$ is invertible and anti-ample.  We shall call such a  surface an \emph{anticanonical surface}. If $d\ge 2$, then the anticanonical morphism
morphism factors through  $\bar S$. For $d\ge 3$ the second factor is an embedding of
$\bar S$ in a projective space of dimension $d$; for $d=3$ this yields a cubic surface and for $d=4$ a complete intersection of two quadrics.  When  $d=2$,
the second factor realizes $\bar S$  as a double cover of a projective plane ramified along a quartic curve with only simple singularities in the sense of Arnol'd
(accounting for the rational double points on $\bar S$).

Adopting the terminology in \cite{dgk}, we say that $S$ is a
\emph{Fano surface} of degree d if $\omega_S^{-1}$ is ample (but beware
that other authors call this a Del Pezzo surface).
If $S$ is given as a projective plane blown up in $9-d$ points, then
it is Fano precisely when the points in question are distinct, no three lie on a line, no six lie on a conic and no eight lie on a cubic which has a singular point at one of them.
This is equivalent to requiring that $S$ contains no $(-2)$-curves.

\emph{From now on we assume that $S$ is not isomorphic to a smooth quadric.}
We denote the canonical class of $S$ by $k\in\pic (S)$ and its orthogonal complement
in $\pic (S)$ by $\pic_0(S)$.
An element $e\in \pic (S)$ is called  an \emph{exceptional class} of $S$ if
$e\cdot e=e\cdot k=-1$. Every exceptional class is representable by
a unique effective divisor.  A \emph{marking} of $S$ is an
ordered $(9-d)$-tuple of exceptional classes $(e_1,\dots ,e_{9-d})$ on $\pic (S)$
with $e_i\cdot e_j=-\delta_{ij}$. Given a marking, there is a unique
class $\ell\in\pic (S)$ characterized by the property that
$3\ell =-k+e_1+\cdots +e_{9-d}$ and $(\ell,e_1,\dots ,e_{9-d})$ will be basis
of $\pic (S)$. The marking is said to be \emph{geometric} if
$S$  can be obtained by $(9-d)$- successive blowups of a projective plane in such a
manner that $e_i$ is  the class of the  total transform $E_i$ of the exceptional
curve of the $i$th blowup. An $\ell$-marking of $S$ consists of merely
giving the class $\ell$. So if $\Lcal$ is a representative line bundle, then $\Lcal$ is
base point free and defines a birational morphism from $S$ to a projective plane
and the anticanonical system on $S$  projects onto a $d$-dimensional linear
system of cubic curves on this plane.

Since we are interested here in the moduli of Fano surfaces, we usually restrict to the case $d\le 4$:  if $S$ is given as a blown up projective plane then four of the $9-d$
points to be blown up can be used to fix a coordinate system, from which
it follows that we have a fine  moduli space  $\Mcal_{m,d}^\circ$ of marked Fano
surfaces of degree $d$ that is isomorphic to an affine open subset of $(\PP^2)^{5-d}$.

\emph{From now on we assume that $d\le 6$.}
With Manin we observe that then the classes
$e_i-e_{i+1}$, $i=1,\dots ,8-d$ and $\ell-e_1-e_2-e_3$
make up a basis of $\pic_0(S)$ and can be thought of as a system of simple roots
of a root system $R_{9-d}$. This root system is of type $E_8$,  $E_7$, $E_6$, $D_5$,
$A_4$ and $A_2+A_1$ respectively. The roots that have fixed inner product with $\ell$
make up a single $\Scal_{9-d}$-orbit and we label them accordingly:
\begin{enumerate}
\item[(0)] $h_{ij}:= e_i-e_j$,  $(i\not= j)$,
\item[(1)] $h_{ijk}:= \ell-e_i-e_j-e_k$ with $i,j,k$ pairwise distinct,
\item[(2)] $(2\ell-e_1-e_2-e_3-e_4-e_5-e_6-e_7)+e_i$, denoted $h_i$ when $d=2$,
For $d=3$, this is only makes sense for $i=7$ and we then may write $h$ instead.
\item[(3)] $-k-e_i$ ($d=1$ only).
\end{enumerate}
Notice that $h_{ij}=-h_{ji}$, but that in $h_{ijk}$ the order of the subscripts is irrelevant.
\[
\DynkinEEEE{h_{123}}{h_{12}}{h_{23}}{h_{34}}{h_{45}}{}{h_{8-d,9-d}}
\]
The marking defined by $(e_1,\dots ,e_{9-d})$ is geometric if and only if
for every $(-2)$-curve $C$ its intersection product with each of the simple roots
is not positive.

The Weyl group $W(R_{9-d})$ is precisely the group of orthogonal transformations of
$\pic (S)$ that fix $k$.  It acts simply transitively of the markings.
In particular it acts on $\Mcal^\circ_{m,d}$ and the quotient variety
$\Mcal_d^\circ:=W(R_{9-d})\bs \Mcal_{m,d}^\circ$ can be interpreted as the coarse moduli space of Fano surfaces of degree $d$. The orbit space of $\Mcal_{m,d}^\circ$ relative to the permutation group of the $e_1,\dots ,e_{9-d}$ (a Weyl subgroup
of type $A_{8-d}$) is the moduli space $\Mcal_{\ell,d}^\circ$ of $\ell$-marked Fano surfaces of degree $d$.

\subsection*{Completion of the moduli spaces by means of GIT}

Fix a $3$-dimensional complex vector space $A$ and a generator $\alpha\in\det (A)$.
We think of $\alpha$ as a translation invariant $3$-vector field on $A$. If $f\in \sym^3 A^*$ is a cubic form on $A$, then the contraction of $\alpha$ with $df$, $\iota_{df}\alpha$, is a $2$-vector field on $A$ that is
invariant under scalar multiplication and hence defines a
$2$-vector field on $\PP(A)$. We thus obtain an isomorphism between $\sym^3 A^*$ and
$H^0(\omega_{\PP (A)}^{-1})$.

Let $d\in \{2,3,4\}$.  A $(d+1)$-dimensional linear quotient  $V$ of
$\sym^3 A$ defines a linear subspace  $V^*\subset \sym^3 A^*$, i.e., a linear
system of cubics on $\PP(A)$ of dimension $d$. If we suppose that
this system does not have a fixed component, then its base locus consists of $9-d$ points (multiplicities counted) and blowing up this base locus produces an $\ell$-marked Del Pezzo surface $S$ with the property that
$H^0(\omega_S^{-1})^*$ can be identified with $V$  (we excluded $d=1$ here because then the  base locus has $9$ points and we get a rational elliptic surface). If we specify an order for the blowing up, then $S$ is even geometrically marked. The quotient surface $\bar S$ obtained from contracting $(-2)$-curves is more canonically defined as it can be described in terms of the rational map $\PP(A)\dashrightarrow \PP (V)$: for $d=3,4$  it is the image of this map and for $d=2$ the Stein factorization realizes $\bar S$ as a double cover $\PP(V)$ ramified over a quartic curve.

The condition that the linear system has no fixed component defines an open subset  $\Omega_d\subset G_{d+1}(\sym^3 A^*)$ whose complement has codimension $>1$. Over
$\Omega_d$ we have a well-defined $\ell$-marked family $\bar \Scal_d/\Omega_d$ to which the $\SL (A)$-action lifts. Any $\ell$-marked anticanonical surface is thus obtained so that we have a bijection between the \emph{set} of isomorphism classes of
$\ell$-marked Del Pezzo surfaces and the \emph{set} of  $\SL(A)$-orbits in
$\Omega_d$. It is unlikely that this can be lifted to the level of varieties and we therefore we invoke  geometric invariant theory. We begin with defining the line bundle that is central to this paper.

\begin{definition}
If $f:\Scal\to B$ is a family of anticanonical surfaces of degree $d$, then its
\emph{determinant bundle} $\Det (\Scal/B)$ is the line bundle over $B$ that is the dual of the determinant of  the rank $9-d$ vector bundle $R^1f_*\omega_{\Scal/B}^{-1}$ (so this assigns to  a Del Pezzo surface $S$, the line $\det H^0(\omega_S^{-1})^*$).
\end{definition}

Thus we have a line bundle $\Det (\Scal_d/\Omega_d)$. Its fiber over
the $(d+1)$-dimensional subspace $V^*\subset \sym^3 A^*$ is
the line $\det (V)$ and hence  the fiber of the ample bundle
$\Ocal_{G_{d+1}(\sym^3 A^*)}(1)$. A  section of
$\Ocal_{G_{d+1}(\sym^3 A^*)}(k)$ determines a section of $\Det^{\otimes k}(\Scal_d/\Omega_d)$ and any section of the latter so occurs because the complement of $\Omega_d$ in $G_{d+1}(\sym^3 A^*)$  has codimension $>1$. Since the action of $\SL (A)$ on $\sym^3 A^*$ is via $\PGL (A)$ (the center $\mu_3$ of $\SL (A)$ acts trivially), we shall regard this as a  representation of the latter.
Consider the subalgebra of $\PGL(A)$-invariants in homogeneous coordinate ring of $G_{d+1}(\sym^3 A^*)$,
\[
R_d^\pt:=\left( \oplus_{k=0}^{\infty} H^0(\Ocal_{G_{d+1}(\sym^3 A^*)}(k))\right)^{\PGL (A)},
\]
The affine cone $\spec (R_d^\pt)$  has the interpretation as the categorical
$\PGL(A)$-quotient of the (affine) Pl\"ucker cone over the
Grassmannian. It may be thought of as  the affine hull of the
moduli space of triples $(S,\ell, \delta)$ with $(S,\ell) $ an $\ell$-marked Del Pezzo surface of degree $d$
and $\delta$ a generator of $\det H^0(\omega_S^{-1})$. Since we shall find that the base of this cone, $\proj (R_d^\pt)$, defines a projective completion of
$\Mcal^\circ_{\ell,d}$, we denote it by $\Mcal^*_{\ell,d}$.
The asserted interpretation of $\Mcal^*_{\ell,d}$ of course requires that
we verify that the orbits defined by Fano surfaces are stable. We will do that in a case by case discussion that
relates this to GIT completions that are obtained in a different manner.
In fact, for each of the three cases $d=2,3,4$ we shall construct
a GIT completion $ \Mcal^*_{d}$ of  $\Mcal^\circ_{d}$ in such a way
that the forgetful morphism $\Mcal^\circ_{\ell,d}\to \Mcal^\circ_{d}$ extends
to a finite morphism of GIT completions $\Mcal^*_{\ell,d}\to \Mcal^*_{d}$.
This description will also help us to identify (and interpret) the boundary strata.

We recall that the proj construction endows $\Mcal^*_{\ell,d}$  for every $k\ge 0$
with a coherent sheaf  $\Ocal_{\Mcal^*_{\ell,d}}(k)$  of rank one whose space of sections is $R_d^k$. We call $\Ocal_{\Mcal^*_{\ell,d}}(1)$ the \emph{determinant sheaf}; it is a line bundle in the orbifold setting.

In what follows, $V_{d+1}$ is a
fixed complex vector space of dimension $d+1$ endowed with  a generator $\mu$ of
$\det (V_{d+1})$. We often regard $\mu$ as a translation invariant
$(d+1)$-polyvector field on $V_{d+1}$.

\subsection*{Degree $4$ surfaces in projective $4$-space}
If a  pencil of quadrics in $\PP(V_5)$ contains a smooth quadric, then
the number of singular members of this pencil (counted with multiplicity) is $5$.
According to Wall \cite{wall:git} the geometric invariant theory for intersections of quadrics is as follows: for a plane $P\subset \sym^2V_5^*$,
$[\wedge^2 P]\in\PP(\wedge^2 (\sym^2V_5^*))$ is $\SL (V_5)$-stable (resp.\   $\SL (V_5)$-semistable)  if and only if the divisor on $\PP(P)$ parameterizing singular members is reduced (resp.\ has all its multiplicities $\le 2$). A semistable pencil belongs to minimal orbit if and only if
its members can be  simultaneously  diagonalized. So a stable pencil is represented by a pair
$\la Z_0^2+Z_1^2+Z_2^2+Z_3^2+Z_4^2,  a_0Z_0^2+a_1Z_1^2+a_2Z_2^2+a_3Z_3^2+a_4Z_4^2\ra$
with $a_0,\dots ,a_4$ distinct. This is equivalent to  the
corresponding surface $S_P$ in $\PP (V_5)$ being smooth. The minimal strictly semistable orbits allow at most two pairs of coefficients to be equal. In case we have only one pair of equal coefficients, $S_P$ has two $A_1$-singularities and in case we have two such pairs, four. The fact  that these singularities come in pairs can be `explained' in terms of the $D_5$-root system in the Picard group of a Del Pezzo surface of degree $4$: a $A_1$-singularity is resolved by a single blowup with a  $(-2)$-curve as exceptional curve whose class is a root in the Picard root system. The roots perpendicular to this root make up a root system of type $D_4+A_1$ and the class of the companion $(-2)$-curve will sit in the $A_1$-summand.  Besides,  a minimal strictly semistable orbit with $2$ resp.\ $4$ $A_1$-singularities  is adjacent  to a semistable orbit without such $A_1$-pairs and represented by  a pair of quadrics one of which is defined by $Z_0^2+Z_1^2+Z_2^2+Z_3^2+Z_4^2$ and the other by
$Z_0Z_1+a_1Z_1^2+a_2Z_2^2+a_3Z_3^2+a_4Z_4^2$ resp.\
$Z_0Z_1+Z_2Z_3+a_3Z_3^2+a_4Z_4^2$.

The center $\mu_5$ of $\SL (V_5)$ acts  acts faithfully by scalars on $\wedge^2(\sym^2V_5^*)$
and for that reason the $\SL (V_5)$-invariant part of the homogeneous coordinate
ring of $\Gr_2(\sym^2V_5^*)$ lives in degrees that are multiples of $5$:
\[
S^\pt_4:= \oplus_{k=0}^\infty S^k_4,\quad S^k_4:= H^0(\Ocal_{\Gr_2(\sym^2V_5^*))}(5k))^{\SL(V_5)}.
\]
We obtain  a projective completion $\Mcal_4^*:=\proj S^\pt_4$  of $\Mcal_4^\circ$ with twisting sheaves  $\Ocal_{\Mcal_4^*}(k)$  such that $S^k_4=H^0(\Ocal_{\Mcal_4^*}(k))$. The singular complete intersections are parameterized by a hypersurface in
$\Gr_2(\sym^2 V_5^*)$. Since the Picard group of this Grassmannian is
generated by $\Ocal_{\Gr_2(\sym^2V_5^*))}(1)$, this discriminant is defined by
a section of some $\Ocal_{\Gr_2(\sym^2V_5^*))}(20)$  and so
$B_4:=\Mcal_4^*-\Mcal_4^\circ$ is defined by a section of $\Ocal_{\Mcal_4^*}(4)$.
\\

Suppose we are given a surface $S\subset \PP(V_5)$ defined
by a pencil of quadrics. So $S$ determines a line $\Phi_S$ in $\wedge^2(\sym^2V^*_5$. Any generator $F_1\wedge F_2\in \Phi_S$ and $u\in V_5^*$  determine
a $2$-vector field on $V$ by  $\iota_{du\wedge dF_1\wedge dF_2}\mu$. This
$2$-vector field is invariant under scalar multiplication and tangent to the cone over $S$.
Hence it defines a $2$-vector field on $S$, or equivalently, an element of $H^0(\omega_S^{-1})$. The map thus defined is an isomorphism
\[
V_5^*\otimes \Phi_S\cong H^0(\omega_S^{-1}).
\]
By taking determinants we get an identification of $\Phi_S^{5}\cong \det H^0(\omega_S^{-1})$. We may think of $\Phi_S^5$ as the quotient of the line  $\Phi_S$ by the center $\mu_5$ of $\SL (V_5)$. Thus $\spec(S^\pt_4)$ may be regarded as the affine hull
of the moduli space of pairs $(S,\delta)$ with $S$ a Del Pezzo surface of degree $4$
and $\delta$ a generator of $\det(H^0(\omega_S^{-1}))$. We conclude:

\begin{proposition}
We have a natural finite embedding $S^\pt_4\subset R^\pt_4$ of graded $\CC$-algebras so that  the forgetful morphism $\Mcal^\circ_{\ell,4}\to \Mcal^\circ_{4}$ extends to a finite morphism of GIT completions $\Mcal^*_{\ell,4}\to \Mcal^*_{4}$ (and the notions of semistability coincide in the two cases)
and $\Ocal_{\Mcal_4^*}(1)$ is the determinant sheaf.
\end{proposition}

\subsection*{Cubic surfaces}
Following Hilbert the  cubic surfaces in $\PP(V_4)$ that are stable (resp.\ semistable) relative to the $\SL (V_4)$-action are those that have an $A_1$-singularity
(resp.\ $A_2$-singularity) at worst. There is only one  strictly
semistable minimal orbit and that is the one  that has three $A_2$-singularities.

The center $\mu_4$ of $\SL (V_4)$ acts faithfully by scalars on $\sym^3V_4^*$
and so the $\SL (V_4)$-invariant part of the homogeneous coordinate
ring of $\sym^3V_4^*$ lives in degrees that are multiples of $4$:
\[
S^\pt_3:= \oplus_{k=0}^\infty S^k_3,\quad S^k_3:= H^0(\Ocal_{\PP(\sym^3V_4^*)}(4k))^{\SL(V_4)}.
\]
We thus find the projective completion $\Mcal_3^*:=\proj S^\pt_3$  of $\Mcal_3^\circ$ with twisting sheaves  $\Ocal_{\Mcal_3^*}(k)$  such that $S^k_3=H^0(\Ocal_{\Mcal_3^*}(k))$. The discriminant hypersurface in the linear system of degree $d$ hypersurfaces in $\PP^n$ has degree
$(n+1)(d-1)^n$. So the singular cubic surfaces are parameterized by a hypersurface in
$\PP(\sym^3 V_4^*)$ of  degree $32$. The stable locus $\Mcal_3^\circ\subset
\Mcal_3\subset \Mcal_3^*$ is the complement of a single point.
Furthermore, $B_3:=\Mcal_3^*-\Mcal_3^\circ$ is defined by a section of
$\Ocal_{\Mcal_3^*}(8)$.
\\

Let $S\subset \PP(V_4)$ be a cubic surface defined by a line $\Phi_S$ in $\sym^3V^*_4$. Proceeding as in the degree $4$ case we find for  that a generator $F\in \Phi_S$
and an $u\in V^*$ the  expression $\iota_{du\wedge dF}\mu$ defines a $2$-vector field
on $S$ and that we thus get an isomorphism
\[
V_4^*\otimes \Phi_S\cong H^0(\omega_S^{-1}).
\]
By taking determinants we get an identification
of $\Phi_S^{4}\cong \det H^0(\omega_S^{-1})$. We think of $\Phi_S^{4}$ as the quotient of the line  $\Phi_S$ by the  center $\mu_4$ of $\SL (V_4)$ and conclude as before:

\begin{proposition}
We have a natural finite embedding $S^\pt_3\subset R^\pt_3$ of graded $\CC$-algebras so that  the forgetful morphism $\Mcal^\circ_{\ell,3}\to \Mcal^\circ_3$ extends to a finite morphism of GIT completions $\Mcal^*_{\ell,3}\to \Mcal^*_{3}$ (and the notions of semistability coincide in the two cases)
and $\Ocal_{\Mcal_3^*}(1)$ is the  determinant sheaf.
\end{proposition}

\subsection*{Quartic curves}
The case of degree $2$ is a bit special because $W(E_7)$
has a nontrivial center (of order two). The center leaves invariant the
(isomorphism type of the) surface: it acts as an involution and only changes
the marking. The latter even disappears if we only remember the fixed point set of
this involution, the quartic curve. Van Geemen \cite{dolgort} observed that
the marking of the Del Pezzo surface then amounts to a principal level two structure
on the quartic curve (this is based on the fact that $W(E_7)$ modulo its center
is isomorphic to the symplectic group $\Sp(6,\ZZ/2)$). Since a smooth quartic curve
is a canonically embedded genus three
curve, $\Mcal_2^\circ$ can also be interpreted  as the moduli space of
nonhyperelliptic  genus three curves with principal level two structure (here we ignore
the orbifold structure).

The projective space  $\PP(\sym^4 V_3^*)$  parameterizes  the quartic
curves in the projective plane $\PP(V_3)$.
The geometric invariant theory relative  its $\SL (V_3)$-action is  as follows:
a quartic curve is stable if and only if it has singularities no worse
than of type $A_2$.
A quartic is unstable if and only if it has a point of multiplicity $\ge 3$
(or equivalently, a $D_4$-singularity or worse) or consists of a cubic plus an
inflectional tangent. The latter gives generically a $A_5$-singularity, but such a
singularity may also appear on a semistable quartic, for instance
on the union of two conics having a point in common where they intersect with
multiplicity $3$. Let us, in order to understand  the incidence relations,  review
(and redo) the classification of nonstable quartics.
\\

A plane quartic curve $C$ that is not stable
has a singularity of type $A_3$ or worse. So it has an equation of the form $cy^2z^2+yzf_2(x,y)+f_4(x,y)$ with $f_2$ and
$f_4$ homogeneous. Consider its orbit under the $\CC^\times$-copy in $\SL (V_3)$
for which $t\in\CC^\times$ sends $(x,y,z)$ to $(x,ty,t^{-1}z)$.
If we let $t\to 0$, then the equation tends to $cy^2z^2+ax^2yz+bx^4$, where
$f_2(x,0)=ax^2$ and $f_4(x,0)=bx^4$. We go through the possibilities.

If $c=0$, then $C$ has a triple point and the equation $ax^2yz+bx^4$ is
easily seen to be unstable. We therefore assume that $c=1$ and we denote the limit
curve by $C_0$.

If $a^2-4b\not=0\not=b$, then $C_0$ is made up of two nonsingular conics meeting
in two distinct points with a common tangent (having therefore a $A_3$-singularity at each) and the original singularity
was of type $A_3$.

If $a\not=0=b$, then we have the same situation except that
one of the conics has now degenerated into a union of two lines.

The most interesting case is when $a^2-4b=0\not=b$. Then
$C_0$ is a double nonsingular conic and in case $C\not= C_0$, $C$ has a singularity of type $A_k$ for some $4\le k\le 7$. The case of an $A_7$-singularity occurs for the curve $C_1$ given by $(yz+x^2)(yz+x^2+y^2)$:  it consist of two nonsingular conics meeting in a single point with multiplicity $4$. This is also the most degenerate
case next after $C_0$: any $\SL (V_3)$-orbit that has $C_0$ is in its closure
is either the orbit of $C_0$ or has $C_1$ in its closure.
So although a double conic does not yield a Del Pezzo surface, the corresponding
point of $\Mcal_2^*$ is uniquely represented by a geometrically marked Del Pezzo surface with a $A_7$-singularity.

On the other hand,
the condition $a=b=0$ (which means that $f_2$ and $f_4$ are divisible by $y$
so that we have a cubic plus an inflectional tangent or worse), gives
the limiting curve defined by $y^2z^2=0$, which is clearly unstable.

We shall later find that the ambiguous behaviour of a $A_5$-singularity
reflects a feature of the $E_7$-root system: this system contains two Weyl group equivalence classes of subsystems
of type $A_5$: one type is always contained in a $A_7$-subsystem (the semistable case)
and the other is not (the unstable case).
\\
Since the center $\mu_3$ of $\SL (V_3)$ acts faithfully by scalars on $\sym^4V_3^*$
we have as algebra of invariants
\[
S^\pt_2:= \oplus_{k=0}^\infty S^k_2,\quad S^k_2:=
H^0(\Ocal_{\PP(\sym^4V_3^*)}(3k))^{\SL(V_3)}.
\]
Thus $\Mcal_2^*:=\proj \CC[\sym^4 V_3^*]^{\SL(V_3)}$
is a projective completion of $\Mcal_2^\circ$. It comes
with twisting sheaves  $\Ocal_{\Mcal_2^*}(k)$  such that
$S^k_2=H^0(\Ocal_{\Mcal_2^*}(k))$.
Let us write $\Mcal_2^\circ\subset \Mcal_2\subset \Mcal_2^*$ for the stable locus;
this can be interpreted as the moduli space of marked Del Pezzo surfaces of degree $2$
with $A_2$-singularities at worst. Its complement in $\Mcal_2^*$ is of dimension one.
Since the singular quartics make up a hypersurface of degree $27$ in $\PP(\sym^4V_3^*)$,
$B_2:=\Mcal_2^*-\Mcal_2^\circ$ is defined by a section of $\Ocal_{\Mcal_2^*}(9)$. In
particular, $B_2$ is a Cartier divisor.
\\

Let $C$ be a quartic curve in $\PP(V_3)$ defined by the line $\Phi_C\subset \sym^4V_3^*$. If $F\in \Phi_C$ is a generator, then a double cover $S$ of
$\PP(V_3)$ totally ramified along $C$ is defined by $w^2=F$ in $V_3\times\CC$ (more precisely, it is $\proj$ of the graded algebra obtained from $\CC[V_3]$ by adjoining to it a root of $F$). Then for every $u\in V_3^*$, the $2$-vector field $w^{-1}\iota_{du\wedge dF}\mu$ defines a section of  $\omega_S^{-1}$.
We thus get an isomorphism $V_3^*\otimes w^{-1}dF\cong H^0(\omega_S^{-1})$.
If we take the determinants of both sides, we find that
$(w^{-1}dF)^3$ determines a generator of  $\det H^0(\omega_S^{-1})$. So
$F^{-3} (dF)^6$ gives one of  $(\det H^0(\omega_S^{-1}))^2$, in other words,
we have a natural isomorphism $\Phi_C^3\cong (\det H^0(\omega_S^{-1}))^2$.
That the square of the determinant appears here reflects the fact that the
central element $-\mathbf{1}$ of $W(E_7)$ induces an involution in $S$
which acts as the scalar $-1$ on $\det  H^0( \omega_{S}^{-1})$. We obtain:

\begin{proposition}\label{prop:detc2system}
We have a natural finite embedding $S^\pt_2\subset R^\pt_2$ of graded $\CC$-algebras so that  the forgetful morphism $\Mcal^\circ_{\ell,2}\to \Mcal^\circ_{2}$ extends to a finite morphism of GIT completions $\Mcal^*_{\ell,2}\to \Mcal^*_{2}$ (and the notions of semistability coincide in the two cases)
and $\Ocal_{\Mcal_2^*}(1)$ is the square of the determinant sheaf.
\end{proposition}

\subsection*{Completion of the moduli space of marked Fano surfaces}
We have produced for $d=4,3,2$, a GIT completion $\Mcal^*_d$ of
$\Mcal^\circ_d$ that we were able to identify with a finite quotient  of
$\Mcal^*_{\ell,d}$. This implies that $\Mcal^*_{\ell,d}$ contains $\Mcal^\circ_{\ell,d}$ as
an open dense subset and proves that every point of $\Mcal^*_{\ell,d}$ can be represented by an $\ell$-marked Fano surface.

We define a completion $\Mcal^*_{m,d}$ of the moduli space $\Mcal^\circ_{m,d}$ of marked Fano surfaces of degree $d$ simply as the normalization of
$\Mcal^*_{d}$ in $\Mcal^\circ_{m,d}$. This comes with an action of $W(R_{9-d})$
and the preceding discussion shows that $\Mcal^*_{\ell,d}$ can be identified
with the orbit space of $\Mcal^\circ_{m,d}$ by the permutation group of the
$e_1,\dots ,e_{9-d}$ (a Weyl subgroup of type  $A_{8-d}$).

\section{Coble's covariants}\label{sect:coblecov}
In this section we assume that the degree $d$ of a Del Pezzo surface is at
most 6 (we later make further restrictions).
\\

Let $(S;e_1,\dots ,e_{9-d})$ be a geometrically marked Del Pezzo surface. Recall that we have a class $\ell\in \pic (S)$ characterized by the property that
$-3\ell+e_1+\cdots +e_{9-d}$ equals the canonical class $k$. Let us choose a line bundle
$\Lcal$ on $S$ which represents $\ell$:  $H^0(\Lcal)$ is then of dimension 3 and if we denote its dual by $A$, then the associated linear system defines a birational morphism $S\to \PP(A)$ which has $E=E_1+\cdots +E_{9-d}$ as its exceptional divisor. The direct image of $\Lcal$ on $\PP(A)$ is still a line bundle (namely $\Ocal_{\PP(A)}(1)$, but we continue to denote  this bundle by $\Lcal$).

We claim that there is a natural identification
\[
\omega^{-1}_S\cong\Lcal^3(-E)\otimes\det A.
\]
To see this, we note that if $p\in \PP(A)$ and $\lambda\subset A$ is the line defined by $p$, then the tangent space of $\PP(A)$ at $p$ appears in the familiar exact sequence
\[
0\to \CC\to \Hom (\lambda ,A)\to T_p\PP(A)\to 0,
\]
from which it follows that $\det T_p\PP(A)=\lambda^{-3}\det (A)$ (we often omit the $\otimes$-symbol when lines or line bundles are involved).
So the anticanonical bundle $\omega_{\PP(A)}^{-1}$ of $\PP(A)$ is naturally identified with
$\Lcal^3\otimes\det (A)$.
Since $S\to \PP(A)$ is the blowup with exceptional divisor $E$, we see that  the above identification makes  $\omega^{-1}_S$ correspond to $\Lcal^3(-E)\otimes\det (A)$.

The following simple lemma will help us to understand  Coble's covariants.

\begin{lemma}
For a Del Pezzo surface $S$, the determinant  lines  of  the vector spaces
$H^0(\Ocal _E\otimes \Lcal^3)\otimes \det (A)$ and $V(S):=H^0(\omega_S^{-1})^*$ are canonically isomorphic.
\end{lemma}
\begin{proof}
The identification $\omega^{-1}_S\cong\Lcal^3(-E)\otimes\det A$ above gives rise to the short exact sequence
\[
0\to \omega_S^{-1} \to\Lcal^3\otimes\det A\to \Ocal_E\otimes \Lcal^3\otimes \det A\to 0.
\]
This yields an exact sequence on $H^0$ because $H^1(\omega_S^{-1})=0$. If we take into account that $H^0(\Lcal^3)=\sym^3H^0(\Lcal)= \sym^3 A^*$, then we find the exact sequence
\[
0 \to  V^* \to \sym^3 A^*\otimes\det A \to H^0(\Ocal_E\otimes \Lcal^3)\otimes \det A\to 0.
\]
Since $\dim (\sym^3 A^*)=10$ and $\det( \sym^3 A^*)=(\det A)^{-10}$, the determinant
of the middle term has a canonical generator. This identifies the determinant of
$V$ with the one of the  right hand side.
\end{proof}

It will be convenient to have a notation for the one dimensional vector space appearing in the preceding lemma: we denote
\[
L(S,E):=\det (H^0(\Ocal _E\otimes \Lcal^3\otimes\det A))=
\det (H^0(\Ocal _E\otimes \Lcal^3))\otimes (\det A)^{9-d},
\]
so that the lemma asserts that $L(S,E)$ may be identified with $\det V(S)$.

We continue with $(S;e_1,\dots ,e_{9-d})$ and $\Lcal$.  If $p_i$ denotes the image point of $E_i$, then the geometric fiber of $\Lcal$ over $p_i$ is $\lambda_i:=H^0(\Lcal\otimes\Ocal_{E_i})^*$ (a one dimensional vector space). So $L(S,E)=
(\lambda_1\cdots \lambda_{9-d})^{-3}\otimes \det (A)^{9-d}$.
For  $e_i,e_j,e_k$ distinct, the map defined by componentwise inclusion
\[
\lambda_i\oplus\lambda_j\oplus\lambda_k\to A
\]
is a linear map between 3-dimensional vector spaces. It is an isomorphism
if $p_i,p_j,p_k$  are not collinear. Hence the corresponding map on the third exterior powers, yields an element
\[
|ijk|\in \lambda_i^{-1}\lambda_j^{-1}\lambda_k^{-1}\det A
\]
that is nonzero in the Fano case (recall that we usually omit the $\otimes$-sign when lines are involved).
Notice that the line
$\lambda_i^{-1}\lambda_j^{-1}\lambda_k^{-1}\det A$ attached to $\Lcal$ only depends on the marked surface $(S;e_1,\dots ,e_{9-d})$ and not on the choice of
the $\Lcal$:  as said above $\Lcal$ is unique up to isomorphism and
the only possible ambiguity therefore originates from the action of $\CC^\times$ in the fibers  of $\Lcal$. But it is clear that this $\CC^\times$-action is trivial on this line.
For  $e_{i_1},\dots , e_{i_6}$ distinct, we also have a  linear map between $6$-dimensional vector spaces
\[
\lambda_{i_1}^2\oplus\cdots \oplus \lambda_{i_6}^2\to \sym^2A.
\]
Since $\det (\sym^2A)=(\det A)^4$, this  defines a determinant
\[
|i_1\cdots i_6|\in \lambda_{i_1}^{-2}\cdots\lambda_{i_6}^{-2} (\det A)^{4}.
\]
It is nonzero if and only if  $p_1,\dots ,p_6$ do not lie on a conic, which is the case
when $S$ is Fano.
The elements $|ijk|$ and $|i_1\cdots i_6|$ just introduced will be referred to as \emph{Coble factors}.

\subsection*{Action of the Weyl group on the Coble factors}
We now assume that $S$ is a Fano surface of degree $\le 6$.
Another marking of $S$ yields another $\ell'$ and hence another line $L(S,E')$. Nevertheless
they are canonically isomorphic to each other since both have been identified with
$\det V$. For what follows it is important to make this isomorphism concrete.
We will do this for the case that the new marking is the image of the former under the
reflection in $h_{123}$.  So $\ell'= 2\ell-e_1-e_2-e_3$, $e'_1=\ell-e_2-e_3$ ($E'_1$ is the strict transform of $\overline{p_2p_3}$), $e'_2$ and $e'_3$ are expressed in a likewise  manner and $e'_i=e_i$ for $i>3$. We represent $\ell'$ by
\[
\Lcal':=\Lcal^2(-E_1-E_2-E_3)
\]
 so that $A'=H^0(\Lcal')^*$. Before proceeding, let us see what happens if we do this twice, that is, if we apply $h_{123}$ once more:
\[
\Lcal''=\Lcal'^2(-E'_1-E'_2-E'_3)=\Lcal^4(-2E_1-2E_2-2E_3-E'_1-E'_2-E'_3).
\]
The line bundle  $\Lcal^3(-2E_1-2E_2-2E_3-E'_1-E'_2-E'_3)$ is trivial (a generator
is given by a section of $\Lcal^3$ whose divisor is  the triangle spanned by
$p_1,p_2,p_3$) and so if $I$ denotes its (one dimensional) space of sections, then
$\Lcal''$ is identified with $\Lcal\otimes I$. We note that for $i>3$, the restriction map
$I\to\lambda_i^{-3}$ is an isomorphism of lines and that we also have a natural isomorphism
$I\to (\lambda_1\lambda_2\lambda_3)^{-1}$, which, after composition with the inverse of $|123|$ yields an isomorphism $I\to \det (A)^{-1}$.

The space of sections of $\Lcal'$ is  the space of quadratic forms on $A$ that are zero on $\lambda_1$, $\lambda_2$ and $\lambda_3$. This leads to an exact sequence
\[
0\to \lambda^2_1\oplus  \lambda^2_2\oplus  \lambda^2_3\to \sym^2 A\to A'\to 0.
\]
The exactness implies  that
\[
\det A'=(\det A)^4(\lambda_1\lambda_2\lambda_3)^{-2}.
\]
We have $(\lambda'_i)^{-1}=H^0(\Lcal^2(-E_1-E_2-E_3)\otimes \Ocal_{E'_i})$
by definition. For $i>3$, this is
just the space of quadratic forms on the line $\lambda_i$, i.e., $\lambda_i^{-2}$ and so
$\lambda'_i=\lambda_i^2$ in that case. For $i=1$,
\[
(\lambda'_1)^{-1}=H^0(\Lcal^2(-E_1-E_2-E_3)\otimes \Ocal_{E'_1})=
H^0(\Lcal^2\otimes \Ocal_{\overline{p_2p_3}})(-(p_2)-(p_3)))
\]
is the space of quadratic forms on
$\lambda_2+\lambda_3$ that vanish on each summand, i.e.,
$\lambda_2^{-1}\lambda_3^{-1}$.  Thus $\lambda'_1=\lambda_2\lambda_3$ and likewise
 $\lambda'_2=\lambda_3\lambda_1$,  $\lambda'_3=\lambda_1\lambda_3$.
Notice that $\lambda''_i$ is naturally identified with $\lambda_i\otimes I^{-1}$. Thus
\begin{multline*}
L(S,E')= (\lambda'_1\cdots \lambda'_{9-d})^{-3} (\det A')^{9-d}=\\
=(\lambda_1\cdots \lambda_{9-d})^{-6}
(\lambda_1\lambda_2\lambda_3)^{-2(9-d)}(\det A)^{4(9-d)}=\\
 =L(S,E)(\lambda_1^{-1}\lambda_2^{-1}\lambda_3^{-1}\det A)^{21-2d}
 \prod_{i=4}^{9-d}(\lambda_i^{-3}\det A).
\end{multline*}
The identifications above of $\lambda_1\lambda_2\lambda_3$, $\lambda_i^3$ ($i>3$)
and $\det (A)$ with $I^{-1}$ show that the twisting line
$(\lambda_1^{-1}\lambda_2^{-1}\lambda_3^{-1}\det A)^{21-2d}
\prod_{i=4}^{9-d}(\lambda_i^{-3}\det A)$ has a canonical generator $\delta$.
This generator can be expressed in terms of our Coble factors as
\[
\delta:=|123|^9\prod_{i=4}^{9-d}(|12i| |23 i| |31i|)
\]
(which indeed lies in
$(\lambda_1^{-1}\lambda_2^{-1}\lambda_3^{-1}\det A)^{21-2d}
\prod_{i=4}^{9-d}(\lambda_i^{-3}\det A)$).

\begin{proposition}\label{prop:comparison}
The isomorphism $L(S;E)\cong L(S',E')$ defined above coincides with the isomorphism that we obtain from the identification of domain and range with $\det V$.
\end{proposition}
\begin{proof}
Choose generators $a_i\in\lambda_i$ and write $x_1,x_2,x_3$ for the basis of $A^*$ dual to $a_1,a_2,a_3$. The basis $(a_1,a_2,a_3)$ of $A$ defines
a generator $a_1\wedge a_2\wedge a_3$ of $\det A$. This determines an  isomorphism
$\phi: \omega_S^{-1}\cong \Lcal^3(-E)$. That isomorphism fits in the exact sequence
\[
0\to V^*\to \sym^3 A^*\to \oplus_{i=1}^{9-d}\lambda_i^{-3}\to 0.
\]
The middle space has the cubic monomials  in $x_1,x_2,x_3$ as a basis. The triple
$(x^3_1,x^3_2,x^3_3)$ defines
a basis dual to  $(a_1^3,a_2^3,a_3^3)\in \lambda_1^{3}\oplus\lambda_2^{3}\oplus\lambda_3^{3}$. It follows that we  have an exact subsequence
\begin{equation}\label{eqn:subseq}
0\to V^*\to K\to \oplus_{i=4}^{9-d}\lambda_i^{-3}\to 0,
\end{equation}
where $K\subset  \sym^3 A^*$ is the span of the cubic monomials that are not a third power.
This yields an identification
\[
(\lambda_4\cdots \lambda_{9-d})^{-3}\cong
\det V \det(\la x_1^2x_2,\dots , x_2x_3^2,x_1x_2x_3\ra ).
\]

We now do the same for $\Lcal'=\Lcal^2(-E_1-E_2-E_3)$. The space $A'$ comes with
a basis $(a'_1=a_2a_3,a'_2=a_3a_1,a'_3=a_1a_2)$ which is dual to the basis $(x_2x_3, x_3x_1, x_1x_2)$ of  $H^0(\Lcal^2(-E_1-E_2-E_3))$. The monomial $x_1x_2x_3$ is the obvious generator of $I=H^0(\Lcal^3(-E'_1-E'_2-E'_3-2E_1-2E_2-2E_3)$ so that have an
associated isomorphism
\[
\phi':=\phi\otimes x_1x_2x_3: \omega_S^{-1}\to (\Lcal')^3(-E').
\]
We fit  this in the exact sequence
\[
0\to H^0(\omega_S^{-1})\to H^0((\Lcal')^3)\to  \oplus_{i=1}^{9-d}\lambda'{}^{-3}_i\to 0.
\]
The  vector space of sections of the middle term has the cubic monomials
in $x_2x_3,x_3x_1, x_1x_2$ as a basis. The lines $\lambda'_1, \lambda'_2,\lambda'_3$ are
spanned by $a_2a_3$, $a_3a_1$ and $a_1a_2$ respectively and $\lambda'_i$ is for $i=4,5,\dots ,9-d$ spanned by $a_i^2$. So $(x_2x_3)^3$ spans $\lambda'{}^{-3}_1$\and similarly for $\lambda'{}^{-3}_2$ and $\lambda'{}^{-3}_3$. The space spanned by cubic monomials in $x_2x_3,x_3x_1, x_1x_2$ that are not pure powers is just
$K':=x_1x_2x_3K$. It follows that  we get an exact subsequence analogous to the sequence (\ref{eqn:subseq}):
\begin{equation}\label{eqn:subseq'}
0\to V^*\to K' \to \oplus_{i=4}^{9-d}\lambda'{}^{-3}_i\to 0,
\end{equation}
where the embedding $V^*\to K'$ is the composite of
$V^*\to K$ and the isomorphism $K\cong K'$ given by multiplication by $x_1x_2x_3$.
This identifies $(a'_i)^{-3}\in\lambda'{}^{-3}_i$ with $(x_1x_2x_3)(a_i).a_i^{-3}\in\lambda_i^{-3}$.
From this we deduce that the generator
$(a'_1\cdots a'_{9-d})^{-3} (a'_1\wedge a'_2\wedge a'_3)^{9-d}$
of $(\lambda'_1\cdots \lambda'_{9-d})^{-3}(\det A')^{9-d}$ corresponds to  $\prod_{i=4}^{9-d}(x_1x_2x_3)(a_i)$ times the generator
$(a_1\cdots a_{9-d})^{-3} (a_1\wedge a_2 \wedge a_3)^{9-d}$ of
 $(\lambda_1\cdots \lambda_{9-d})^{-3}(\det A)^{9-d}$.
Since $\delta=\prod_{i=4}^{9-d}(x_1x_2x_3)(a_i)$, this proves that the isomorphism $L(S;E)\cong \det V\cong L(S,E')$ sends the generator of the former to  $\delta$ times the generator of the latter.
\end{proof}

We next determine how the Coble factors for $\ell'$ are expressed in terms of those of
$\ell$. We retain our $0\not= a_i\in\lambda_i$ and write $x_1,x_2,x_3$
for the basis of $A^*$ dual to $a_1,a_2,a_3$ as before.
We identify $|ijk|$ resp. \ $|i_1\cdots i_6|$ with the $a_1\wedge a_2\wedge a_3$-coefficient of $a_i\wedge a_j\wedge a_k$ resp.\ the
$a_1^2\wedge a^2_2\wedge a^2_3\wedge a_2a_3\wedge a_3a_1\wedge a_1a_2$-coefficient of $a_{i_1}^2\wedge\cdots \wedge a^2_{i_6}$. Here are some typical cases, where it is assumed that the free indices are distinct and $>3$:
\begin{align*}
|123| &=1,\\
|12k| &=\la x_3\, |\, a_k\ra ,\\
|1JFK|  &=\la x_2\wedge x_3\, |\,  a_j\wedge a_k\ra , \\
|ijk | &=\la  x_1\wedge x_2\wedge x_3\, |\,  a_i\wedge a_j\wedge a_k\ra , \\
|123ijk|&=
\la x_2x_3\wedge x_3x_1\wedge x_1x_2\, |\,  a^2_i\wedge a^2_j\wedge a^2_k\ra,\\
|12ijkl|&=
\la x_3^2\wedge x_2x_3\wedge x_3x_1\wedge  x_1x_2\, |\,  a^2_i\wedge a^2_j\wedge a^2_k\wedge a_l^2\ra,\\
|1ijklm|&=
\la x_2^2\wedge x_3^2\wedge x_2x_3\wedge x_3x_1\wedge  x_1x_2\, |\,
a^2_i\wedge a^2_j\wedge a^2_k\wedge a^2_l\wedge a_m^2\ra .
\end{align*}
The corresponding expressions for the new marking are converted into
the old marking by the substitutions
\[
a'_i=
\begin{cases}
a_2a_3 &\text{when $i=1$,}\\
a_3a_1 &\text{when $i=2$,}\\
a_1a_2 &\text{when $i=3$,}\\
a_i^2 &\text{when $i>3$.}
\end{cases},\quad
x'_i=
\begin{cases}
x_2x_3 &\text{when $i=1$,}\\
x_3x_1 &\text{when $i=2$,}\\
x_1x_2 &\text{when $i=3$,}\\
x_i^2 &\text{when $i>3$.}
\end{cases}.
\]
We thus find:
\begin{align*}
|123|' &=1=|123|,\\
|12k|' &=x_1x_2(a_k)=|23k| |31k|,\\
|1jk |'  &=-x_1(a_j)x_1(a_k)|1jk|,\\
|ijk|' &= |123ijk|,\\
|123ijk|'&=x_1x_2x_3(a_i). x_1x_2x_3(a_j).x_1x_2x_3(a_k).|ijk|.\\
|12ijkl|'&=x_1x_2(a_i). x_1x_2(a_j).x_1x_2(a_k).x_1x_2(a_l).|12ijkl|,
\end{align*}
The expression for $|1ijklm|'$ does not appear to have a pleasant form: we find  that\begin{multline*}
|1ijklm|'= x_1(a_i)x_1(a_j)x_1(a_k)x_1(a_l)x_1(a_m).\\
.\la x_3^2x_1\wedge x_1x_3^2\wedge x_1x_2x_3\wedge x_2^2x_3\wedge  x_2x_3^2\, \, |\,
a_i^3\wedge a_j^3\wedge a_k^3\wedge a_l^3\wedge a_m^3\ra .
\end{multline*}

\subsection*{The covariants} Here is the definition.

\begin{definition}\label{def:coblevector}
Let $(S; e_1,\dots ,e_{9-d})$ be a  marked  Fano surface
of degree $d\le 6$. A \emph{Coble covariant} is an element of $L(S,E)$ that is a product of Coble factors $|ijk|$ and $|i_1\cdots i_6|$ in such a manner that every unordered pair in $\{ 1,2,\dots ,9-d\}$ appears in  one of these factors.
\end{definition}

This notion also makes sense for a marked Del Pezzo surface, and indeed,
in case the $E_1,\dots , E_{9-d}$ are irreducible (or equivalently, $p_1,\dots ,p_{9-d}$
are distinct), then we adopt this as a  definition. But when this is not the case, this is not
the `right' definition (see Remark \ref{rem:coblevector}).

It is easily verified that Coble covariants exist only when $2\le d\le 5$. In these cases they are as follows:
\begin{description}
\item[$d=5$] There is only one Coble covariant, namely
$|123| |234| |341| |412|$. It is nonzero if and only if no three points are collinear, that is,
if $S$ is a Fano surface.
\item[$d=4$] A typical Coble covariant is $|123| |234| |345| |451| |512|$.  It depends on a cyclic ordering of $\{1,2,\dots ,5\}$, with  the opposite cycle giving the same element. So the number of Coble covariants up to sign is equal to $4!/2=12$. A Coble covariant
can be nonzero even if $S$ has $(-2)$-curves.  For instance, if $(p_1,p_2,p_4)$ and $(p_2,p_3,p_5)$ are collinear but are otherwise
generic then the given Coble covariant is nonzero and $S$ has a $2A_1$-configuration
(i.e., two disjoint $(-2)$-curves).
\item[$d=3$]  We have two  typical cases: one is $|134| |234| |356| |456| |512|
|612|$ and another is $|123| |456| |123456|$. The former type amounts to  dividing the 6-element set $\{e_1,\dots ,e_6\}$ in three equal parts (of two) and cyclically order the three parts (there are $30$ such) and the latter to splitting of $\{e_1,\dots ,e_6\}$ into two equal parts (there are $10$ of these). So there are $40$ Coble covariants up to sign.
\item[$d=2$]  We have two typical cases: $|351| |461| |342| |562| |547| |217| |367|$ (of which there are 30) and $|123456| |127| |347| |567|$ (105 in number).
So up to sign we find 135 cases.
\end{description}

\begin{proposition}
If $S$ is Fano, then the collection of Coble covariants, when considered as elements of  $\det V(S)$, is independent of the marking.
\end{proposition}
\begin{proof}
It is enough to show that the collection is invariant under the reflection in
$h_{123}$.
So in view of Proposition  \ref{prop:comparison} we need to verify that if we make the above substitutions for a Coble covariant relative to $(e'_1,\dots ,e'_{9-d})$, then we get $\delta$ times a Coble covariant
relative to $(e_1,\dots ,e_{9-d})$. This is a straightforward check. We do a few examples.
For $d=3$, we find
\begin{multline*}
|134|'  |234|' |356|' |456|'  |512|' |612|'=\\
=-x_1x_3(a_4). x_2x_3(a_4). x_3(a_5)x_3(a_6)|356|.|123456|.x_1x_2(a_5).x_1x_2(a_6)=\\
=-\delta |124| |356| |123456|,
\end{multline*}
which is indeed the image of $|124| |356| |123456|$. Similarly,
\[
|123|'  |456|' |123456|'=
|123|.|123456|.\delta |456|,
\]
which is the image of $|123| |456| |123456|$.

The other cases are similar and are left to the reader to verify.
\end{proof}

\begin{remark}\label{rem:coblevector}
The notion of a Coble covariant extends to the case of a geometrically marked
Del Pezzo surface. There is not much of an issue here as long as the points
$p_1,\dots ,p_{9-d}$ remain distinct, but when two coalesce  the situation
becomes a bit delicate, since we wish to land in $\det V(S)$. It is clear that if
 $p_2$ approaches $p_1$, then  any Coble covariant involving these points
such as $|123|\in \lambda_1^{-1}\lambda_2^{-1}\lambda_3^{-1}\det A$ tends to $0$.
But we should regard  $|123|$ as an element of
$\det (H^0(\Ocal_{E_1+E_2+E_3}\otimes\Lcal^3))\otimes\det A$ and when $p_2$ tends to $p_1$, then $E_1$ becomes decomposable and of the form $F+E_2$.  The component $F$ is  the strict transform of the exceptional curve of the first blowup and hence a $(-2)$-curve. Let $q$ be the point where $F$ and $E_2$ meet. This corresponds to tangent direction at $p_1$, or equivalently, to a plane $P_q\subset  A$ that contains $\lambda_1$. If $\lambda_2$ moves in $W$ towards $\lambda_1$, then
$H^0(\Ocal_{E_1+E_2})$ becomes $H^0(\Ocal_{F+2E_2})$. The exact sequence
\[
0\to \Ical_{F+E_2}/\Ical_{F+2E_2}\to \Ocal_{F+2E_2}\to \Ocal_{F+E_2}\to 0
\]
induces an exact sequence on sections. Notice that the first term is a constant sheaf on $E_2$. Its fiber over $q$ is $T^*_qF\otimes T^*_qE_2$. This fiber is apparently also the determinant of  $H^0( \Ocal_{F+2E_2})$. Since
$T_qF=\Hom(\lambda, W_q/\lambda)$ and $T_qE_2=\Hom(W_q/\lambda ,A/W_q)$, we
have  $T^*_qF\otimes T^*_qE_2=(\det A)^{-1}\det P_q\otimes \lambda_1$. It follows that
\[
\det (H^0(\Ocal_{F+2E_2+E_3}\otimes\Lcal^3))\otimes\det A
=\lambda_1^{-1}\lambda_3^{-1}\det P_q.
\]
It is in this line where $|123|$ should take its value. (This also explains why in the
next section we need to divide by a discriminant $\Delta (t_q,\dots ,t_{9-d})$.)
\end{remark}

Since every point of $\Mcal_{m,d}^*$ is representable by a marked Del Pezzo surface,
a Coble covariant  can be regarded as a section of the determinant sheaf
$\Ocal_{\Mcal_{m,d}^*}(1)$. It follows from Section \ref{sect:coblecov} that $W(R_{9-d})$ permutes
these sections transitively.

\begin{definition}
The \emph{Coble space}  $\coble_d$ is the subspace of
$H^0(\Ocal_{\Mcal_{m,d}^*}(1))$ spanned by the Coble covariants.
\end{definition}

We shall prove that  for $d=3,4$, $\Ccal_d$ is complete, i.e., all of
$H^0(\Ocal_{\Mcal_{m,d}^*}(1))$,
but we do not know whether that is true when $d=2$.

\begin{remark}
We shall see that $\coble_d$ is an irreducible representation of
$W(R_{9-d})$ of dimension $6, 10, 15$ for resp.\ $d=4,3,2$.
It follows from the discussion in Section \ref{sect:gitmoduli} that any
$W(R_{9-d})$-invariant polynomial of degree $k$ in
the Coble covariants has an interpretation in terms of classical invariant
theory: for $d=4$ we get a $\SL(5)$-invariant of degree $5k$
for pencils of quadrics, for $d=3$, a $\SL(4)$-invariant of degree $4k$
for cubic forms and for $d=2$, a $\SL(3)$-invariant of degree $3k/2$
for quartic forms.

Here is an example that illustrate this. There is only one irreducible
representation of $W(E_{6})$ of degree $10$. This  representation is real and
has therefore a nonzero $W(E_{6})$-invariant quadratic form. According to the
preceding this produces a
$\SL(4)$-invariant of degree $8$ for cubic forms. This is indeed the lowest
degree of a such an invariant.
\end{remark}

\section{Anticanonical divisors with a cusp}\label{sect:cusp}

\subsection*{Anticanonical cuspidal cubics on Del Pezzo surfaces}
Let $S$ be a Del Pezzo surface of degree $d$ that is not isomorphic to a smooth quadric.
Assume  that is also given a reduced anticanonical curve $K$ on $S$ isomorphic to a cuspidal cubic
(notice that if $d=1$ such a curve will not always exist). The curve $K$ is given by a hyperplane  $V_K\subset V=H^0(\omega^{-1}_S)^*$. We write $l_K$ for the line $V/V_K$ so that $l_K^*$ is a line in
$V^*=H^0(\omega_S^{-1})=\Hom (\omega_S,\Ocal_S)$.
The image of $l_K^*$ is $\Hom (\omega_S,\Ocal_S(-K))\subset \Hom (\omega_S,\Ocal_S)$
and hence $l_K$ may be identified with  $H^0(\omega_S(K))$.
So a nonzero $\kappa\in l_K$ can be understood as a rational
$2$-form $\kappa$ on $S$ whose divisor is $K$.  The residue $\res_K(\kappa)$
of  $\kappa$ on the smooth part of $K$ identifies $\pic^o(K)$ with $\CC$ as an algebraic group: we may represent an element  of $\pic^o(K)$ by a difference $(q)-(p)$ and
then $\res_K(\kappa)$ assigns to this element the integral of $\res_K(\kappa)$ along
any arc in $K$ from $p$ to $q$. This identifies $l_K$ with $\Hom(\pic^o(K),\CC)$ or equivalently, $\pic^o(K)$ with $l^*_K$.

Recall that $\pic_0(S)\subset \pic (S)$ denotes the orthogonal complement of the class of $\omega_S^{-1}$. It is then clear that restriction defines a homomorphism
$r:\pic_0(S)\to\pic^o(K)\to l^*_K(\subset V^*)$. We extend $r$ to an algebra homomorphism
\[
r: \sym^\pt \pic_0(S)\to \sym^\pt l^*_K (\subset \sym^\pt V^*).
\]
For $\kappa\in l_K$, we compose this map with the evaluation in $\kappa$ and obtain
an algebra homomorphism
\[
r_\kappa: \sym^\pt \pic_0(S)\to \CC.
\]
Suppose now $S$ geometrically marked by  $(e_1,\dots ,e_{9-d})$ as before.
With the notation of the previous section, we have a line bundle $\Lcal$ on $S$ and an associated contraction morphism $S\to \PP(A)$, where $A=H^0(\Lcal)^*$, with $E_i$ mapping to a singleton. The cuspidal curve $K$ meets $E_i$ in a single point $p_i$
(with intersection number one). It is mapped isomorphically to its image in $\PP(A)$.

The Zariski tangent space of $K$ at its cusp is a line (with multiplicity two, but that
will be irrelevant here). Let $u\in A^*$ be such that $u=0$ defines the corresponding
line in $\PP(A)$. We may then extend $u$ to a coordinate system $(u,v,w)$ for
$A$ such that $K$ is given by the equation $u^2w-v^3$ (this makes $[1:0:0]$ the
unique flex point of $K$ and $w=0$ its tangent line). This coordinate system is  for a
given $u$ almost unique: if $(u,v',w')$ is any other such coordinate system, then
$v'=cv$ and $w'=c^3w$ for some $c\in\CC^\times$.

However, a choice  of a generator $\kappa\in l_K$ singles out a natural choice
$(v,w)$ by requiring that the residue of $\kappa$ on $K$ is the restriction of $d(v/u)$. We put $t:=v/u$ (we should write $v_\kappa,w_\kappa, t_\kappa$, but we
do not want to overburden the notation, let us just remember that $v,w, t$ are homogeneous of degree $1,3,1$ in $\kappa$). The dependence of $v$ and $w$
on $u$ is clearly homogeneous of degree $1$.
The smooth part $K$ is then parameterized by $t\in\CC\mapsto p(t):= [1:t:t^3]$ such that
$dt$ corresponds to $\res_K\kappa$ and $r_\kappa$ sends  $(p(t))-(p(t'))\in \pic^o(K)$ to
$t-t'\in \CC$.

Assume for the moment that the $p_i$'s are distinct (so that the $E_i$'s are irreducible
and the $\lambda_i$'s are distinct). Let us denote the restriction of $u\in A^*$ to $\lambda_i$ by $u_i$. This is clearly a coordinate for $\lambda_i$ and hence a generator of $\lambda_i^{-1}$.
We thus obtain the generator
\begin{align*}
\eps'_\kappa:=(u_1\cdots u_{9-d})^{3} (du\wedge dv\wedge dw)^{-(9-d)}\\ \in \lambda_1^{-3}\cdots\lambda_{9-d}^{-3} (\det A)^{9-d}=L(S,E).
\end{align*}
Since $v$ and $w$ are homogeneous of degree one  in $u$, it follows
that $\eps'_\kappa$ is independent of the choice of $u$. But they are  homogeneous of degree $1$ and $3$ respectively  in $\kappa$, and so $\eps'_\kappa$ is homogeneous of degree $-4(9-d)$ in $\kappa$.

Let us now see which linear forms we get on
$\lambda^3_1\oplus\cdots\oplus\lambda^3_{9-d}$ by restriction of  cubic monomials in
$u,v,w$. They will  be of the
form $(t_1^ku_1^3,\cdots ,t_{9-d}^ku_{9-d}^3)$ for some $k\ge 0$:
for the monomial $u^a v^b w^c$ we have $k=b+3c$. We thus get all integers $0\le k\le 9$ except $8$ (and $3$ occurs twice since $u^2v$ and $w^3$ yield the same restriction):
\[
\begin{array}{ccccc}
u^3\mapsto 1 & u^2v\mapsto t & uv^2\mapsto t^2 &u^2w, v^3\mapsto t^3&uvw\mapsto t^4\\
v^2w\mapsto t^5 &uw^2\mapsto t^6 & vw^2\mapsto t^7 &w^3\mapsto t^9.
\end{array}
\]
If we select $9-d$ such monomials and compute the  determinant of their restrictions
to $\lambda^3_1\oplus\cdots\oplus\lambda^3_{9-d}$, we see that that it is either zero or
equal to $(u_1\cdots u_{9-d})^{3}\det ((t_i)^{k_j})_{1\le i,j\le 9-d}$ for some
$0\le k_1<\cdots <k_{9-d}$. The latter expression lies in $\ZZ[t_1,\dots ,t_{9-d}]$ and is divisible by the discriminant that we get by taking  $9-d$ monomials corresponding to
$1,t,\dots ,t^{8-d}$, namely
\[
\Delta(t_1,\dots ,t_{9-d}):=\prod_{1\le i<j\le 9-d} (t_i-t_j).
\]
In other words, if we regard $t_1,\dots ,t_{9-d}$ as variables, then the $(9-d)$th exterior power over $\ZZ[t_1,\dots ,t_{9-d}]$ of the homomorphism
\[
\ZZ[t_1,\dots ,t_{9-d}]\otimes \sym^3A^*\to \oplus_{i=1}^{9-d}\ZZ[t_1,\dots ,t_{9-d}]u_i^3
\]
has its image generated by $\Delta(t_1,\dots ,t_{9-d}) u_1^3\cdots u_{9-d}^{3}$. That
a division by $\Delta(t_1,\dots ,t_{9-d})$ is appropriate is suggested by Remark
\ref{rem:coblevector} and so we use
\[
\eps_\kappa:=\Delta(t_1,\dots ,t_{9-d})\eps'_\kappa
\]
as a generator of
$\lambda_1^{-3}\cdots\lambda_{9-d}^{-3} (\det A)^{9-d}$ instead. We may rephrase this more sensibly in terms of our exact sequence
\[
0\to V^*\to \sym^3(A^*)(\det A)^{10}\to H^0(\Lcal^3\otimes\Ocal_E\otimes\det A)\to 0.
\]
We see that our coordinates define a basis of the middle term in such a manner that
they make the sequence split:  $9-d$ cubic monomials that yield
$1,t,\dots ,t^{8-d}$ define a partial basis of $\sym^3(A^*)(\det A)^{10}$ whose
$(9-d)$th exterior power maps onto $\eps_\kappa$. Notice that this remains true
if some of the points $p_i$ coalesce, that is, if $S$ is just a geometrically marked
Del Pezzo surface---this is in contrast to $\eps'_\kappa$. Since $\eps_\kappa$ is homogeneous  in $\kappa$ of degree $-4(9-d)+\binom{9-d}{2}=\tfrac{1}{2}d(9-d)$, we have constructed an isomorphism
\[
\eps: l_K^{d(9-d)/2}\cong L(S,E).
\]

\subsection*{Anticanonical cuspidal cubics on Fano surfaces}Let us return to the Fano case.
If $p_i=p(t_i)=[1:t_i:t_i^3]$, then the generator of $\lambda_i$ dual to $u_i$ is clearly $\tilde p_i=(1,t_i,t_i^3)$.
It is not hard to verify that for $i,j,k$ distinct  in $\{1,\dots ,9-d\}$,
we have the following identity in $(\lambda_i\lambda_j\lambda_k)^{-1}\det A$:
\begin{multline*}
\tilde p_i\wedge \tilde p_j\wedge\tilde p_k
=\Delta (t_i,t_j,t_k) (-t_i-t_j-t_k)  (du\wedge dv\wedge dw)^{-1}=\\
=\Delta (t_i,t_j,t_k) r_\kappa(h_{ijk})(du\wedge dv\wedge dw)^{-1}.
\end{multline*}
Similarly we find  for $i_1,\dots ,i_6$ distinct in $\{1,\dots ,9-d\}$ the following identity in
the line $\lambda_{i_1}^{-2}\cdots \lambda_{i_6}^{-2}\det (\sym^2A)$:
\begin{multline*}
\tilde p_{i_1}^2\wedge\cdots \wedge\tilde p_{i_6}^2=
\mp(t_{i_1}+\cdots +t_{i_6})\Delta(t_{i_1},\dots ,t_{i_6})(du\wedge dv\wedge dw)^{-4}=\\
=\pm \Delta(t_{i_1},\dots ,t_{i_6}) r_\kappa(2\ell -e_{i_1}-\cdots -e_{i_6}) (du\wedge dv\wedge dw)^{-4}.
\end{multline*}
In this manner we get for example when $d=3$:
\begin{multline*}
\pm |123456| |123| |456|=\\
=\Delta(t_1,\dots ,t_{6}) r_\kappa(\Delta(R'))u_1^3\cdots u_6^3(du\wedge dv\wedge dw)^{-6}=r_\kappa(\Delta')\eps_\kappa,
\end{multline*}
where $\Delta'$ is the discriminant of the root subsystem of type $3A_2$ given by
$\la h_{12},h_{23}, h_{45},h_{56},h_{123}, h\ra$, and similarly
\[
\pm |134| |234| |356| |456| |512| |612|=
r_\kappa(\Delta'')\eps_\kappa ,
\]
where $\Delta''$ is the $3A_2$-discriminant of
$\la h_{12}, h_{134},h_{34}, h_{356}, h_{56}, h_{125}\ra$. In this manner  to each Coble covariant there is associated the discriminant of a
$A_2^3$-subsystem of the $E_6$ root system and vice versa (there are indeed
$40$ such subsystems). In similar fashion
we find that for $d=2$ a Coble covariant is the discriminant of a
$A_1^7$-subsystem of the $E_7$ root system and vice versa
(there are $135$ such); the two typical cases yield
\begin{align*}
\pm |135| |146| |234| |256| |457| |127| |367| &=
r_\kappa(h_{135}h_{146}h_{234}h_{256}h_{457}h_{127}h_{367})
\eps_\kappa,\\
\pm |123456| |127| |347| |567| &=
r_\kappa(h_{12}h_{34}h_{56}h_{127}h_{347}h_{567}h_7)
\eps_\kappa.
\end{align*}
For $d=4$, we do not get the discriminant of a root subsystem. The best way to describe this case is perhaps by just giving a typical case in terms
of the standard representation of $D_5$ in $\RR^5$ as in Bourbaki (\cite{bourbaki:lie}), where the roots
are $\pm\eps_i\pm\eps_j$ with $1\le i<j\le 5$. One such case is
\[
\prod_{i \in \ZZ/5}(\eps_i-\eps_{i+1})(\eps_i+\eps_{i+1})=
\prod_{i \in \ZZ/5} (\eps_i^2-\eps_{i+1}^2).
\]
There are indeed $12$ such elements up to sign.

Finally, we observe that for $d=5$, we get
\[
\pm |123||234||341||412|= r_\kappa(\Delta)\eps_\kappa,
\]
where $\Delta$ is the discriminant  of the full $A_4$-system $\la h_{12}, h_{23}, h_{34}, h_{123}\ra$.

This makes it clear that in all these cases $\eps: l_K^{d(9-d)/2}\cong L(S,E)$ is independent of the marking once we identify $L(S,E)$ with $\det V$. Notice that the polynomials defining Coble covariants indeed have the predicted degree
$\frac{1}{2}d(9-d)$.
\\

We can restate this as follows. Consider the Lobatschevki lattice
$\Lambda_{1,9-d}$ whose basis elements are denoted $(\ell, e_1,\dots ,e_{9-d})$
(the inner product matrix is in diagonal form with $(+1,-1,\dots ,-1)$ on the diagonal).
Put $k:=-3\ell +e_1+\cdots +e_{9-d}$ and let
$(h_{123}=3\ell-e_1-e_2-e_3, h_{12}=e_1-e_2,\dots , h_{8-d,9-d}=e_{8-d}-e_{9-d})$.
This is a basis of $k^\perp$ that is at the same time a
root basis of a root system $R_{9-d}$.

If $R$ is a root system, then let us denote its root lattice $Q(R)$,  by
$W(R)$ its Weyl group, by
$\hfrak(R):=\Hom (Q(R) ,\CC)$ the Cartan algebra on which $R$ is defined.
Let  $\hfrak(R)^\circ\subset  \hfrak(R)$ stand for the reflection hyperplane complement (which, in the parlance of Lie theory, is the set of its \emph{regular elements}). We abbreviate the  projectivizations of these last two spaces by $\PP(R)$ and $\PP(R)^\circ$. In the presence of a nondegenerate $W(R)$-invariant symmetric
bilinear form on $Q(R)$ we tacitly identify $\hfrak(R)$ with its dual.

So $Q(R_{9-d})=k^\perp$. It is clear that a marking of a Del Pezzo surface amounts to an isomorphism $\pic (S)\cong \Lambda_{1,9-d}$ which sends the canonical class to $k$ (and hence $\pic_0(S)$ to $Q(R_{9-d})$). These  isomorphisms are simply transitively permuted by the Weyl group $W(R_{9-d})$.
If we are given a marked Fano surface $(S;e_1,\dots ,e_{9-d})$ of degree $d$ and a rational $2$-form $\kappa$ on $S$ whose divisor a cuspidal curve $K$, then we can associate to these data an element of $\hfrak(R_{9-d})$ by
\[
Q(R_{9-d})\cong \pic_0(S)\to \pic^o(K) {\buildrel {r_\kappa}\over \longrightarrow} \CC.
\]
It is known that we land in $\hfrak(R_{9-d})^\circ$ and that we thus obtain a bijection between the set of isomorphism classes of systems $(S;e_1,\dots ,e_{9-d}; \kappa)$ and the points of $\hfrak(R_{9-d})^\circ$. This isomorphism is evidently homogeneous of degree one: replacing  $\kappa$ by $c\kappa$ multiplies the image by a factor $c$. In other words, the set of isomorphism classes of systems $(S;e_1,\dots ,e_{9-d}; K)$, where  $(S;e_1,\dots ,e_{9-d})$ is a marked Fano surface of degree $d$ and $K$ is a cuspidal anticanonical curve on $K$ can be identified with
$\PP(\hfrak_{9-d})^\circ$ (where we have abbreviated  $\hfrak(R_{9-d})$ by
$\hfrak_{9-d}$) in a such a manner that
if $l$ is a line in $\hfrak_{9-d}$ associated to  $(S;e_1,\dots ,e_{9-d}; K)$, then $l$ gets identified with  the line $H^0(\omega_S(K))$.

We sum up the preceding in terms of the forgetful  morphism
$\PP(\hfrak_{9-d})^\circ\to \Mcal_{m,d}^\circ$:

\begin{theorem}\label{thm:projquotient}
Assume that $d\in\{ 2,3,4,5\}$. Then the  forgetful  morphism from $\PP(\hfrak_{9-d})^\circ$ to
$\Mcal_{m,d}^\circ$  is surjective and flat. It is covered by a natural isomorphism between the pull-back of the determinant bundle and $\Ocal_{\PP(\hfrak_{9-d})^\circ}(\frac{1}{2}d(9-d))$ and under this isomorphism the Coble covariants form a single $W(R_{9-d})$-orbit, which up to a constant common scalar factor is as follows:
\begin{enumerate}
\item[(d=5)] the discriminant of the $A_4$-system (a polynomial of degree $10$),
\item[(d=4)] the $W(D_5)$-orbit of the degree 10
polynomial $\prod_{i \in \ZZ/5} (\eps_i^2-\eps_{i+1}^2)$,
\item[(d=3)] the discriminants of subsystems of type $3A_2$ (a $W(E_6)$-orbit of
polynomials of degree $9$),
\item[(d=2)] the discriminants of subsystems of
type $7A_1$ (a $W(E_7)$-orbit of polynomials of degree $7$).
\end{enumerate}
In particular, the Coble space $\coble_d$ can be identified with the linear span of the
above orbit of polynomials.
\end{theorem}

The Weyl group representation $\coble_d$ was, at least  for $d=2$ and $d=3$,
already considered by Coble \cite{coblebook}, although the notion of a Coxeter group
was not available to him. We shall consider these representations in more
detail in Section \ref{sect:coblerep}. In that section we also investigate  the
rational map  $\PP(\hfrak_{9-d})\dashrightarrow \Mcal^*_{m,d}$ defined
by the morphism $\PP(\hfrak_{9-d})^\circ\to \Mcal^\circ_{m,d}$.
\\

Theorem \ref{thm:projquotient} presents (for $d=2,3,4,5$) the moduli space $\Mcal_{m,d}^\circ$ as a flat $W(R_{9-d})$-equivariant quotient of $\PP(R_{9-d})^\circ$ with $(d-2)$-dimensional fibers. Admittedly that description is somewhat indirect from a  geometric point of view.  We will here offer in the next two subsections a somewhat more concrete characterisation when  fiber and base are positive dimensional (so for $d=3$ and
$d=4$). A fiber  is then irreducible and we show that the $(d-1)$
$W(R_{9-d})$-invariant vector fields of lowest degree span in $\PP(R_{9-d})^\circ$
a $(d-2)$-dimensional foliation whose leaves are the fibers of $\PP(R_{9-d})^\circ\to
\Mcal_{m,d}^\circ$.  We first do the case $d=3$.

\subsection*{The universal parabolic curve of a cubic surface}
We begin with recalling a classical definition.

\begin{definition}
Let $S$ be a Fano surface  of degree $3$.
A point of $S$  is said to be \emph{parabolic} if there is a cuspidal anticanonical curve
on $S$ that has its cusp singularity at  that point.
\end{definition}

We may think of the surface $S$ as anticanonically embedded in $\PP^3$ as a smooth cubic surface.
If $F(Z_0,Z_1,Z_2,Z_3)$ is a defining equation, then the locus  of parabolic points is precisely the part of $S$ where $S$ meets its  \emph{Hessian
surface} (defined by $\det (\partial^2F/\partial Z_i\partial Z_j)$) transversally. So this is a nonsingular curve on $S$ that need not be closed (in fact, it isn't: a
point  of $S$ where its tangent plane intersects $S$ in a union of a conic and a line tangent lies in the Zariski boundary of the parabolic curve). If we fix a marking  for $S$, so that is determined a point of $p\in\Mcal^\circ_{m,3}$, then
the parabolic locus can be identified with
the fiber of $\PP(R_{9-d})^\circ$ over $p$.

We now return to the situation of the beginning of this section, where we essentially have a fixed cuspidal
cubic curve $K$ in the projective plane $\PP^2$ whose smooth part has in terms of affine coordinates has the parameter form $(v,w)=(t, t^3)$ (this puts the cusp at infinity and the unique flex point at the origin). The points $p_1,\dots ,p_6\in\PP^2_\ell$ that we blow up in order to produce $S$ lie on $K_\reg$; we denote their
$t(=v)$-coordinates $t_1,\dots ,t_6$. The system $(S;e_1,\dots ,e_6;K)$ defines
a point $\tilde p\in \PP(E_6)^\circ$ and $(t_1,\dots ,t_6)$ describes a point of
$\hfrak(E_6)^\circ$ that lies over $\tilde p$.
The vector fields $X$ on $\PP^2$ with the property
that $X$ be tangent to $K$ at the points $p_1,\dots ,p_6$ make
up a vector space of dimension two. It contains the field $X_E=v\partial /\partial v+3w\partial /\partial w$, which is
tangent to $K_\reg$ everywhere (it generates a $\CC^\times$-action on $\PP^2$ that preserves $K$). If $X$ is in this vector space,  then
\[
\hat X:=\sum_i X(p_i)\frac{\partial}{\partial t_i}
\]
is a tangent vector of $\hfrak(E_6)$ at $(t_1,\dots ,t_6)\in \hfrak(E_6)^\circ$.
For $X_E$ this yields $\hat X_E=\sum_{i=1}^6 t_i\frac{\partial}{\partial t_i}$,
in other words, we get the Euler field of $\hfrak(E_6)$ at $(t_1,\dots ,t_6)$.

\begin{lemma}
If $X$ is not proportional to $X_E$, then the line in $\PP(E_6)$ through
$\tilde p$ that is defined by $\hat X$
is tangent to the fiber of  $\PP(E_6)^\circ\to \Mcal_{m,3}^\circ$ at $\tilde p$.
\end{lemma}
\begin{proof}
This is mostly a matter of geometric  interpretation.  If we view $\hat X=(X(p_1),\dots ,X(p_6))$  as an  infinitesimal displacement of the point configuration
$(p_1,\dots ,p_6)$ in $\PP^2$,  then $\hat X$ does not effectively deform the corresponding Fano surface, because $X$ is an infinitesimal automorphism of $\PP^2$.
But if we view  $\hat X$ as an  infinitesimal displacement $(p_1,\dots ,p_6)$ then it will
induce a nontrivial line field (a priori with singularities) unless $X$ is  tangent to $K$.
This last condition is equivalent to $X$ being proportional to $X_E$.
\end{proof}

We now calculate the resulting field on $\PP(E_6)^\circ$.
A vector field $X$ on $\PP^2$ has in our affine coordinates $(v,w)$ the form
\[
X=(a_0+a_1v+a_2w +c_1v^2+c_2vw)\frac{\partial}{\partial v}+
(b_0+b_1v+b_2w +c_1vw+c_2w^2)\frac{\partial}{\partial w}
\]
Since we may calculate modulo $X_E$ we assume $b_2=0$. The condition that $X$ be
tangent to $C$ at $p_i$ amounts to:
$3t^2(a_0+a_1t+a_2t^3+ c_1t^2+c_2t^4)=b_0+b_1t+c_1t^4+c_2t^6$
for $t=t_i$, or equivalently, that
$2c_2t^6+ 3a_2t^5+2c_1t^4+3a_1t^3+3a_0t^2-b_1t-b_0$ has
(the distinct) zeroes $t=t_i$ for $i=1,\dots ,6$. This means that
\begin{align*}
\frac{a_2}{c_2}=-\tfrac{2}{3}\sigma_1,\;
\frac{c_1}{c_2}=\sigma_2,\;
\frac{a_1}{c_2}=-\tfrac{2}{3}\sigma_3,\;
\frac{a_0}{c_2}=\tfrac{2}{3}\sigma_4,\;
\frac{b_1}{c_2}=2\sigma_5,\;
\frac{b_0}{c_2}=-2\sigma_6,\;
\end{align*}
where $\sigma_i$ stands for the $i$th symmetric function of $t_1,\dots ,t_6$.
So we may normalize $X$ by taking $c_2=1$. The value of $X$ in $p_i$ is
in terms of the $t$-coodinate its $x$-component and hence equal to
\[
(a_0+a_1t_i+c_1t_i^2+a_2t_i^3+c_2t_i^4)\frac{\partial}{\partial t}=
(\tfrac{2}{3}\sigma_4-\tfrac{2}{3}\sigma_3 t_i +\sigma_2 t_i^2
-\tfrac{2}{3}\sigma_1 t_i^3+ t_i^4)\frac{\partial}{\partial t}.
\]
It follows that
\begin{equation*} \tag{$*$}
\hat X=\sum_{i=1}^6(\tfrac{2}{3}\sigma_4-\tfrac{2}{3}\sigma_3 t_i +\sigma_2 t_i^2
-\tfrac{2}{3}\sigma_1 t_i^3+ t_i^4)\frac{\partial}{\partial t_i}.
\end{equation*}

We see in particular that if regard $\hat X$ as a vector field ($(t_1,\dots ,t_6)$ varies), then it is homogeneous of degree $3$.
The space of homogenenous vector fields of degree $3$ on $\hfrak (E_6)$
is as a $W_6$-representation space isomorphic to
$\sym^4(\hfrak^* (E_6))\otimes \hfrak (E_6)$; this has a one-dimensional
space of invariants and contains no other $W(E_6)$-invariant one-dimensional
subspace. It follows that $\hat X$ is $W(E_6)$-invariant and is characterized by this
property up to a constant factor.

\begin{corollary}
The fibration $\PP(E_6)^\circ \to\Mcal_{m,3}^\circ$ integrates the one
dimensional foliation defined by a $W(E_6)$-invariant vector field that is
homogeneous of degree three.
\end{corollary}

A natural way to produce  such an invariant vector field is to take the nonzero
$W(E_6)$ invariants  polynomials
$f_2,f_5$ on $\hfrak(E_6)$ of degree $2$ and $5$ (these are unique up to a constant
factor); since $f_2$ is nondegenerate we can choose coordinates
$z_1,\dots ,z_6$ such that $f_2=\sum_i z_i^2$. The gradient vector field relative to $f_2$,
\[
\nabla f_5=\sum_{i=1}^6 \frac{\partial f_5}{\partial z_i}
\frac{\partial}{\partial z_i}
\]
is a $W(E_6)$-invariant homogeneous vector field of degree $3$. So
we can restate the preceding corollary as

\begin{theorem}
Let $\hfrak$ denote the natural representation space of a Coxeter group
$W$ of type $E_6$. The  natural $W$-invariant rational dimension-one
foliation on $\PP(\hfrak)^\circ$ of degree $3$ (i.e., the one defined by the gradient of a
nonzero invariant quintic form  with  respect a
nonzero (hence nondegenerate) invariant quadratic form on $\hfrak$) is
algebraically integrable and has a leaf space that is in a $W$-equivariant manner
isomorphic to the moduli space of marked smooth cubic surfaces.
\end{theorem}

\subsection*{Moduli of degree $4$ Del Pezzo surfaces}
This is a  slight modification of the argument for the degree $3$ case.
We have one point less and so by letting $t_6$ move over the affine line
we may regard formula ($*$) as defining a one parameter family of
vector fields on $\hfrak(D_5)$. For $i\ge 1$, we have $\sigma_i(t_1,\dots ,t_6)=
\sigma_i(t_1,\dots ,t_5)+t_6\sigma_{i-1}(t_1,\dots ,t_5$) and so we immediately
see that we this is a linear family spanned by the two vector fields
\begin{align*}
\hat X_3=&\sum_{i=1}^5(\tfrac{2}{3}\sigma_4-\tfrac{2}{3}\sigma_3 t_i +\sigma_2 t_i^2
-\tfrac{2}{3}\sigma_1 t_i^3+ t_i^4)\frac{\partial}{\partial t_i},\\
\hat X_2=&\sum_{i=1}^5(\tfrac{2}{3}\sigma_3-\tfrac{2}{3}\sigma_2t_i +\sigma_1 t_i^2
-\tfrac{2}{3} t_i^3)\frac{\partial}{\partial t_i}.
\end{align*}
The subscript indicated of course the degree. Since $\hat X$ is $W(E_6)$-invariant,
$\hat X_3$ and $\hat X_2$ will be invariant under the $W(E_6)$-stabilizer of $e_6$, that is, $W(D_5)$. The $W(D_5)$-invariant vector fields on $\hfrak (D_5)$ form a free module on the polynomial algebra of  $W(D_5)$-invariant functions.
The latter algebra has its generators $f_2,f_4,f_5,f_6,f_8$ in degrees indicated by the subscript. The generator $f_2$ is a nondegenerate quadratic form and the module of invariant vector fields is freely generated by  the gradients of the $f_i$ relative to $f_2$,
$X_{i-2}:=\nabla f_i$ (these have the  degree indicated by the subscript; $X_0$ is the Euler field).  We conclude that the plane distribution on $\PP(D_5)^\circ$ spanned by the vector fields $X_2$ and $X_3$ is also defined by $\hat X_2$ and $\hat X_3$.
We conclude :

\begin{theorem}
Let $\hfrak$ denote the natural representation space of a Coxeter group
$W$ of type $D_5$. The  natural $W$-invariant rational dimension-two
foliation on $\PP(\hfrak)^\circ$ defined by  the gradients of a
nonzero invariant forms of degree $4$ and $5$  with  respect a
nonzero (hence nondegenerate) invariant quadratic form on $\hfrak$ is
algebraically integrable and has a leaf space that is in a $W$-equivariant manner
isomorphic to the moduli space of marked Fano surfaces of degree $4$.
\end{theorem}

\begin{remark}
The Frobenius integrability  is remarkable, because
it tells us that the degree four vector field $[X_2,X_3]$ does not involve the  degree
four generator $X_4$. It is of course even more remarkable that it is algebraically so (in the sense that its leaves are the fibers  of a morphism). It makes one wonder how often this happens.
For instance one can ask: given a Coxeter arrangement complement $\hfrak^\circ$ and a positive integer $k$, when is the distribution on $\hfrak^\circ$ spanned by the subset of homogeneous invariant generating vector fields  whose degree is  $\le k$
algebraically  integrable?
\end{remark}

\section{Coble's representations}\label{sect:coblerep}

This section discusses the main properties of the representations of a Weyl group
of type $D_5$, $E_6$  or $E_7$ that we encountered in Theorem \ref{thm:projquotient}.

\subsection*{Macdonald's irreducibility theorem}
We will use the following beautiful (and easily proved!) theorem of
MacDonald \cite{macdonald} which states that the type of representation
under consideration is irreducible.

\begin{proposition}[Macdonald]\label{prop:macdonald}
Let $R$ be a root system, $\hfrak$ the complex vector space it spans and
$S\subset R$ a reduced root subsystem. Then the $W(R)$-subrepresentation
of $\sym^{|S|/2} \hfrak$ generated by the discriminant of $S$ is irreducible.
In particular, the Coble
representations of type $E_6$ and
$E_7$ are irreducible.\qed
\end{proposition}
\begin{proof}
Since the proof is short, we reproduce it here.
If $L\subset \sym^{|S|/2} \hfrak$ denotes the line spanned by the discriminant of
$S$, then  $W(S)$ acts on $L$ with the sign character. In fact, $L$ is the entire eigensubspace of $\sym^{|S|/2} \hfrak$ defined by that character is, for if $G\in \sym^{|S|/2}\hfrak$ is such that $s(G)=-G$ for every reflection $s$ in $W(S)$, then $G$ is zero on
each reflection hyperplane of $W(S)$ and hence divisible by the discriminant
of $S$. Since $G$ and the discriminant have the same degree,
$G$ must be proportional to it.

Let $V=\CC[W(R)]L$ be the  $W(R)$-subrepresentation of $\sym^{|S|/2} \hfrak$ generated by $L$. We must prove that  every $W$-equivariant map $\phi:V\to V$ is given by a scalar. From the preceding it follows that $\phi$ preserves $L$ and so is given on $L$ as multiplication by a scalar,  $\lambda$ say. Then $\phi-\lambda\mathbf{1}_V$ is zero on $L$ and hence zero on $V$.
\end{proof}

Unfortunately Macdonald's theorem does not come with an effective way
to compute the degree of such representations  and that is one
of several good reasons to have a closer look
at them. (The Coble representations had indeed been  considered by Coble
and presumably by others before him. Their irreducibility and their degrees
were known at the time.)

It will be convenient (and of course quite relevant for the application we have in mind)
to work  for $d=2,3,4$ with the Manin model of the $R_{9-d}$ root system as sitting
in the Lobatchevski lattice $\Lambda_{1,9-d}$
so that $\hfrak_{9-d}:=\hfrak(R_{9-d})$ is the orthogonal element of $k=-3\ell+e_1+\cdots +e_{9-d}$.
For $d=2,3,4$ we have a corresponding (Coble) representation $\coble_{d}=
\coble (R_{9-d})$, which in case $d=3$ resp.\ $d=2$ is  spanned by the
discriminants subsystems of type $3A_2$ resp.\ $7A_1$. We may regard  $\coble_{d}$ as a linear system of hypersurfaces  of degree $10$ ($d=4$), $9$ ($d=3$) or $7$ ($d=2$) in $\PP(\hfrak_{9-d})$. Among our goals is to compute  the dimension of this  system and to investigate
its separating  properties.

The Manin basis recognizes one particular weight, namely the orthogonal
projection  of $\ell$ in $\hfrak_{9-d}$. Its $W(R_{9-d})$-stabilizer is the symmetric group
$\Scal_{9-d}$ of $e_1,\dots ,e_{9-d}$ (a Weyl subgroup of type $A_{8-d}$).
We shall denote by $\pi_{9-d}: \CC\otimes\Lambda_{1,{9-d}}\to \hfrak_{9-d}$ the orthogonal projection. So $\pi_{9-d}(e_i)=e_i+\frac{1}{3}k$.

\subsection*{The Coble representation of a Weyl group of type $E_6$}
Here $R=R_6$  and the representation of $W(E_6)$ in question is the subspace
$\coble_3\subset\sym^9\hfrak_6$ spanned by the discriminants of subsystems of
type $3A_2$ of $R$. Following Proposition  \ref{prop:macdonald},
this representation is irreducible. We shall prove that $\dim \coble_3=10$ and that
quotients of elements of $\coble_3$ separate the isomorphism types of
cubic surfaces.

\begin{lemma}\label{lemma:e6properties}
The Weyl group $W(R)$ acts transitively on the collection of ordered triples of
mutually orthogonal roots.  If  $(\alpha_1,\alpha_2,\alpha_3)$ is such
a triple, then
\begin{enumerate}
\item[(i)] there is a root $\alpha\in R$ perpendicular
to each $\alpha_i$ and this root is unique up to sign,
\item[(ii)] the roots $\alpha_1$, $\alpha_2$, $\alpha_3$, $\alpha$
belong to (unique) subsystem of type $D_4$,
\item[(iii)] there are precisely two subsystems of type $3A_2$ containing
$\{\alpha_1,\alpha_2,\alpha_3\}$ and these two subsystems are interchanged by $s_\alpha$.
\end{enumerate}
\end{lemma}
\begin{proof}
The  transitivity assertion and the properties (i) and (ii) are known and a proof goes like this: the orthogonal complement of a root $\alpha_1\in R$  is a subsystem $R'\subset R$ of type $A_5$, the orthogonal complement of a root $\alpha_2\in R'$ in $R'$  is a subsystem $R''\subset R'$ of type $A_3$ and the orthogonal complement of a root $\alpha_3\in R''$ is a
subsystem $R'''\subset R''$ of type $A_1$. Since all the root systems encountered have the property that their Weyl group acts transitively on the roots, the first assertion follows.
Notice that we proved (i) at the same time. The remaining properties now only need
to be verified for a particular choice of $(\alpha_1,\alpha_2,\alpha_3)$.

We take this triple to be  $(h_{12},h_{34},h_{56})$. Then
we may take $\alpha =h$ and we see that these are roots of
the $D_4$-system spanned by $h_{12}$, $h_{34}$, $h_{56}$, $h_{135}$.
The two $3A_2$-subsystems containing $\{h_{12},h_{34},h_{56}\}$ are then easily seen to be
$\la h_{12}, h_{134}\ra \perp \la h_{34}, h_{356}\ra \perp \la h_{56}, h_{125}\ra$ and
$\la h_{12}, h_{156}\ra\perp \la h_{34}, h_{123}\ra
\perp \la h_{56}, h_{345}\ra$.
We observe that $s_h$ interchanges them.
\end{proof}

The following notion is the root system analogue  of  its namesake introduced by
Allcock and Freitag \cite{af}.

\begin{lemmadef}\label{lemdef:cross}
Let $R$ be a root system of type $E_6$ and $S\subset R$ a subsystem of type $3A_2$. If $\alpha\in R$  is not orthogonal to any summand of $S$, then the roots in $S$ orthogonal to $\alpha$ make up a subsystem of type $3A_1$ (which then must meet every summand of $S$). This sets up a bijection between the antipodal pairs $\{ \pm\alpha\}$ that are not orthogonal to any summand of $S$ and $3A_1$-subsystems of $S$.

For $(S,\alpha)$ as above and $S^+$ a set of positive roots for $S$, the
degree nine polynomial $(1-s_\alpha)\Delta (S^+)$ is called a \emph{cross} of $R$.
\end{lemmadef}
\begin{proof}
If we are given a $A_2$-subsystem of $R$, then any root
not in that subsystem is orthogonal to some root in that subsystem. This implies
that in the above definition we can find in each of the three $A_2$-summands a root
that is orthogonal to $\alpha$. Since $\alpha$ is not orthogonal to any summand, this root is unique up to sign and so the roots in $S$ fixed by $s$  form a $3A_1$-subsystem as asserted.

Conversely, if $R'\subset S$ is a subsystem of type $3A_1$, and
$\alpha\in R-R'$ is as  in Lemma \ref{lemma:e6properties}, then $s=s_\alpha$ has
the desired property.
\end{proof}

\begin{lemma}\label{lemma:quintic}
Let $S^+\subset R$  and  $\alpha\in R-S$ be as in Lemma \ref{lemdef:cross}. If $\alpha_1,\alpha_2,\alpha_3$ are the roots in $S^+$ perpendicular to $\alpha$, then the cross
$(1-s_\alpha)\Delta (S^+)$ is divisible by $\alpha_1\alpha_2\alpha_3\alpha$:
\[
(1-s_\alpha)\Delta (S^+)=\alpha_1\alpha_2\alpha_3\alpha  F_1,
\]
and the quotient $F_1\in \sym^5\hfrak_6$  is invariant under the Weyl group of the
$D_4$-subsystem that contains $\alpha_1,\alpha_2,\alpha_3,\alpha$.
\end{lemma}
\begin{proof}
It is clear that both $\Delta (S^+)$ and $s_\alpha\Delta (S^+)=\Delta (s_\alpha S^+)$ are divisible by $\alpha_1\alpha_2\alpha_3$.  It is also clear that $(1-s_\alpha)\Delta (S^+)$
is divisible by $\alpha$.  So $F_1$ is defined as an element of $\sym^5\hfrak$.

We will now prove that there exists a $g\in \GL(\QQ\otimes\Lambda_{1,6})$
which centralizes the Weyl group in question and is such that the transform
of $F_1$ under $g^{-1}$ is a $W(R)$-invariant in $\sym^5\hfrak$. This will clearly suffice.

We may, in view of Lemma \ref{lemma:e6properties},  assume without loss of generality that $(\alpha_1,\alpha_2,\alpha_3,\alpha)=(h_{12},h_{34},h_{56},h)$ so that the $D_4$-subsystem containing these roots is
$\la h_{12}, h_{34}, h_{56}, h_{135}\ra $. Denote by $\hfrak'$ the subspace of $\hfrak$ spanned by these roots.

 We first recall a remarkable result
due to Naruki. The set of exceptional
classes that have inner product $1$ with $\alpha$ is  $\{e_1,\dots ,e_6\}$
(this set and its $s_\alpha$-transform make up what is classically known as a
\emph{double six}). Consider the element
$(1-s_\alpha)\prod_{i=1}^6 \pi_*(e_i)\in \sym^6 \hfrak$. It is clearly
divisible by $\alpha$ and the quotient $F\in \sym^5\hfrak$  will evidently be  invariant under a
Weyl subgroup of $W(R)$ of type $A_5+A_1$. But according to
Naruki (\cite{hunt}, p.\ 235) $F$ is even invariant under all of $W(R)$.

The orthogonal complement of $\hfrak'$ in $\CC \otimes\Lambda_{1,6}$ is spanned by the
members of the `anticanonical triangle'
\[
(\eps_0:=\ell-e_1-e_2,\eps_1:=\ell-e_3-e_4,\eps_2:=\ell-e_5-e_6)
\]
and the intersection $\hfrak\cap \hfrak'{}^\perp$ is spanned by the differences $\eps_0-\eps_1$ and
$\eps_1-\eps_2$.
Let $g\in \GL(\QQ\otimes\Lambda_{1,6})$ be the transformation that is the identity on
$\hfrak'$ and takes $\eps_i$ to $\eps_i+2\eps_{i+1}$ for $i\in\ZZ/3$. This transformation
preserves $\hfrak'{}^\perp$ and hence commutes with all the transformations that preserve $\hfrak'$
and act as the identity on $\hfrak'{}^\perp$. We also note that $g(k)=3k$, that $g$
preserves the orthogonal complement of $k$ and hence $g$ commutes with $\pi$.

We claim that $g\pi(e_1)=-h_{134}\in R$.
One easily checks that $2e_1+\eps_0+k\in \hfrak'$ and so
\begin{multline*}
g(e_1)=g(\half(2e_1+\eps_0+k))-\half g(\eps_0)-\half g(k)=\\
=\half(2e_1+\eps_0+k)-\half(\eps_0+2\eps_{1})-\tfrac{3}{2}k=\\
=e_1-\eps_1-k=-\ell+e_1+e_3+e_4-k=-h_{134}+k.
\end{multline*}
Applying $\pi$ to this identity yields $g\pi(e_1)=-h_{134}$.

We get similar formulas for the $g\pi(e_i)$ and thus find that
\[
g\pi \{e_1,\dots ,e_6\} =\la -h_{134}, -h_{234}\ra\perp
\la-h_{356}, -h_{456}\ra \perp \la-h_{125},-h_{126}\ra.
\]
Notice that the union of this set with $\{ h_{12},h_{34},h_{56}\}$ is a
system of positive roots of a $3A_2$-system.
This union will be our $S^+$. So $g_*$ takes the polynomial
$h_{12}h_{34}h_{56}\prod_{i=1}^6 \pi (e_i)$ to $\Delta (S^+)$. Since $s_\alpha$ commutes with $g$ we have
\begin{multline*}
(1-s_{\alpha})\Delta (S^+)
=h_{12}h_{34}h_{56}g\Big((1-s_{\alpha})  \prod_{i=1}^6 \pi (e_i)\Big)=\\
=h_{12}h_{34}h_{56}g(\alpha F)=h_{12}h_{34}h_{56}\alpha g(F)
\end{multline*}
so that $F_1=g(F)$. This proves the lemma.
\end{proof}

\begin{corollary}\label{cor:crossa1}
For any three pairwise perpendicular roots in $R$ there exists a cross that is
divisible by their product. This cross is unique up to sign and is
also divisible by a root perpendicular to these three. This yields a bijection
between $4A_1$-subsystems of $R$ and antipodal pairs of crosses.\qed
\end{corollary}

\begin{corollary}\label{cor:crossd4}
Let $R'\subset R$  be a subsystem of type $D_4$. The $4A_1$-subsystems of  $R'$  define three crosses up to sign whose sum (up to sign) is zero. These crosses span a $W(R')$-invariant plane in $\coble_3$. A quotient of the discriminants of two
$4A_1$-subsystems of  $R'$ is a quotient of two crosses.\qed
\end{corollary}

We now fix a $D_5$-subsystem $R_o\subset R$. It has precisely $5$ subsystems
of type $D_4$. As we just observed, each of these defines a plane in $\coble_3$. Therefore, the crosses associated to the $4A_1$-subsystems of $R_o$ span a subspace of $\coble_3$ of dimension at most $10$.

\begin{lemma}\label{lemma:d5trans}
For every subsystem $S\subset R$ of type $3A_2$, $S\cap R_o$ is of type
$2A_1+ A_2$ and hence contains three subsystems of type $3A_1$;
for every such $3A_1$-subsystem the associated $4A_1$-subsystem of $R$ is
in fact contained in $R_o$.
Moreover $S\mapsto S\cap R_o$ defines a bijection
between the $3A_2$-subsystems of $R$ and the $2A_1+A_2$-subsystems
of $R_o$ and  $W(R_o)$ acts transitively on both sets.
\end{lemma}
\begin{proof}
It is easy (and left to the reader) to find one subsystem $S\subset R$
of type $3A_2$ such that $R_o\cap S$ has the stated properties.
It therefore suffices to prove the transitivity property.
This involves a simple count: The $W(R_o)$-stabilizer of $R_o\cap S$
of $R_o$ contains $W(R_o\cap S)$ as a subgroup of index two (there is
an element in the stabilizer that interchanges the $A_1$-summands and
is minus the identity on the $A_2$-summand) and so the number of systems
$W(R)$-equivalent to
$R_o\cap S$  is $|W(D_5)|/2|W(A_2+ 2A_1)|=40$.
That is just as many as there are $3A_2$-subsystems of $R$.
\end{proof}

We continue with the $D_5$-subsystem $R_o\subset R$ that we fixed above. Let  $S\subset R$ be any subsystem of type $3A_2$. By Lemma \ref{lemma:d5trans} $R_o\cap S$ is of type  $2A_1+ A_2$. Let be $s_1$, $s_2$, $s_3$ be the three
reflections in the Weyl group of the $A_2$-summand. The two $A_1$-summands
and the antipodal root pair attached to $s_i$ make up
a $3A_1$-subsystem $R_1^{(i)}$ of $R_o\cap S$. Each of these
subsystems is contained in a unique  $4A_1$-subsystem.
Let $s^{(i)}$ denote the reflection in the extra $A_1$-summand.
According to Lemma \ref{lemma:d5trans}, $s^{(i)}\in W(R_o)$, so that
$R^{(i)}_2:= s^{(i)} S$ has the property that $R^{(i)}\cap R_o=R\cap R_o$.

\begin{lemma}\label{lemma:crosssum}
The discriminant $\Delta (S^+)$ is fixed under $s'+s''+s'''$, in other words,
\[
2\Delta (S^+)=(1-s')\Delta (S^+)+(1-s'')\Delta (S^+)+(1-s''')\Delta (S^+),
\]
where we note that the right hand side is a sum of three crosses attached to
subsystems of $R_o$ of type $4A_1$. In particular, $\coble_3$ is generated by
the crosses.
\end{lemma}
\begin{proof}
It is clear that $s_i\Delta (S^+)=-\Delta (S^+)$.
Since $s_3=s_1s_2s_1^{-1}$, we have $R_1'''=s_1(R_1'')$. This implies that
$s'''=s_1s''s_1$ and so
\[
(s''+s''')\Delta (S^+)=(s''+s_1s''s_1)\Delta (S^+)= (1-s_1)s''\Delta (S^+)=(1-s_1)\Delta(s''S^+).
\]
Since $s_1\notin W(s''S)$, the right hand side is a cross.
We claim that this cross equals the cross $(1-s')\Delta (S^+)$ up to sign.
For this it suffices to show that there exist four perpendicular roots such that
each is divisible by three of them. It is clear that $(1-s_1)s''\Delta (S^+)$
is divisible by a root attached to $s_1$ and by the roots in the two
$A_1$-summands of $S^+$ (for these are unaffected by $s''$ and $s_1$). On the other hand,
$(1-s')\Delta (S^+)$ is divisible by a root attached to $s'$  and the roots in the two
$A_1$-summands of $S$. It remains to observe that the roots attached to $s_1$ and
$s'$ are perpendicular.

Thus $(s''+s''')\Delta (S^+)= \pm (1-s')\Delta (S^+)$. Suppose the minus sign holds, so that
$1-s'+s''+s'''$ kills $\Delta (S^+)$. The cyclic permutation $1+s'+s''-s'''$ then
also kills $\Delta (S^+)$ and hence so will $1+s''$. In other words, $\Delta (S^+)$ will be anti-invariant
under $s''$. Since $s''\notin W(S)$, this is a contradiction. Hence the plus sign holds
and the lemma follows.
\end{proof}

\begin{theorem}
The planes  defined by the five subsystems of $R_o$ of type $D_4$
make up a direct sum decomposition of $\coble_3$. In particular $\coble_{3}$ is
the irreducible representation of the $E_6$-Weyl group of degree 10.
\end{theorem}
\begin{proof}
Lemma  \ref{lemma:crosssum} shows that  $\coble_3$ is spanned by the crosses
attached  to $4A_1$-subsystems of $R_o$. So the five planes in question span $\coble_3$
and $\dim\coble_3\le 10$.
The irreducible representations of $W(R)$ of degree $<10$ are the trivial
representation, the sign representation (which are both of degree $1$) and the defining
representation (of degree $6$) and $\coble_3$ is clearly neither of these.
The theorem follows.
\end{proof}

We now determine the common zero set of the Coble covariants. We first make some general remarks that also apply to the $D_5$ and the $E_7$-case.
The zero set of a Coble invariant  is a union of reflection hyperplanes and hence
each irreducible component  of their common intersection, $Z_r\subset \hfrak_r$, is an intersection of reflection hyperplanes. As $Z_r$ is invariant under the Weyl group, so is the collection of its  irreducible components. So an irreducible component is always the translate  of common zero set of a subset of the given root basis of $R_r$. (This subset need not be unique.)

\begin{proposition}\label{prop:z6}
The common zero set $Z_6\subset \hfrak$ of the members of  $\coble_3$
is the union of the linear subspaces that are pointwise fixed by a
Weyl subgroup of type $A_3$.
\end{proposition}
\begin{proof}
We first verify that for any $A_3$-subsystem of $R$, the subspace of $\hfrak$ perpendicular to it is in the common zero set of the members of $\coble_3$.
Since an $A_3$-subsystem is contained in a $D_5$-subsystem, it is in view of
Lemma \ref{lemma:d5trans} enough to show that a $A_2+2A_1$-subsystem
and a $A_3$-subsystem in a $D_5$-system always meet. This is easily verified.

We next show that $Z_6$ is not larger. Any subsystem generated by fundamental roots
that does not contain a $A_3$-system is contained in a subsystem of type $2A_2+A_1$. There is a single Weyl group equivalence class of such subsystems and so
it suffices to  give two subsystems of $R$,  of type $3A_2$ and of type $2A_2+A_1$ that are disjoint. We take $\la h_{12}, h_{23}, h_{45},h_{56},h_{123},h\ra$
and  $\la h_{16}, h_{125}, h_{34}, h_{136}, h_{25}\ra$.
\end{proof}

\begin{question}
Is $\coble_3$ the space of degree $9$ polynomials on $\hfrak$ that vanish on $Z_6$?
This is probably equivalent to the completeness of $\coble_3$ as a linear system on
$\Mcal^*_{m,3}$ (which is known, though in a rather indirect manner, see Remark
\ref{rem:modularforms}).
\end{question}

\subsection*{The Coble representation of a Weyl group of type $E_7$}
The Weyl group $W(E_7)$ decomposes as  $W_+(E_7)\times\{1,c\}$, where
$W_+(E_7)\subset W(E_7)$ is the subgroup of elements that have determinant one
in the Coxeter representation and $c\in W(E_7)$ is minus the identity in
the Coxeter representation. This implies that every irreducible representation of $W(E_7)$ is obtained as an irreducible representation of $W_+(E_7)$ plus a decree as to whether $c$ acts as $1$ or as $-1$.

We know that the representation of $W(E_7)$ defined by
$\coble_2$  (which we recall, is spanned by  products of seven pairwise perpendicular roots of the $E_7$ root system) is irreducible and we want to prove:

\begin{proposition}\label{prop:coble7}
The representation $\coble_2$ of $W(E_7)$ is of degree $15$ and the nontrivial central element of $W(E_7)$ acts as $-1$.
\end{proposition}

It is known that there is just one isomorphism types of irreducible representations of
$W_+(E_7)$ in degree $15$  and so Proposition \ref{prop:coble7}  identifies the isomorphism type of the representation.

In what follows $R$ stands for the root system $R_7$ of type $E_7$.

\begin{lemma}\label{lemma:7a1}
The Weyl group $W(R)$ acts transitively on the collection of $7A_1$-subsystems of
$R$.  If  we are given a subsystem $R'$ of type $2A_1$, then  the roots perpendicular to $R'$ make up a subsystem  of type $A_1+ D_4$. In particular, there is a unique subsystem of type $3A_1$ that contains $R'$ and is orthogonal to a subsystem of type
$D_4$. Conversely, the roots perpendicular to a given subsystem of $R$ of type $D_4$
make up a system of type $3A_1$.
\end{lemma}
\begin{proof}
This lemma is known and the proof is standard. The first assertion follows from the fact that the roots orthogonal to a given root of $R$ form a subsystem of type $D_6$ and the roots orthogonal to a root of a root system of type $D_6$ form a subsystem of type
$D_4+A_1$.

Any root subsystem of $R$ of type $D_4$ is saturated  and so a root basis of this subsystem extends to a root basis of
$R$. As the group $W(R)$ acts transitively on the set of root bases, it also acts transitively on the set of subsystems of type $D_4$.
\end{proof}

So if we have a subsystem $R_1\subset R$ of type $7A_1$, then any two
summands of $R_1$ (making up a subsystem $R'\subset R_1$ of type $2A_1$)
determine a third summand and the remaining $4$ summands will lie in a
$D_4$-subsystem.
In this way we can construct a $2$-dimensional simplicial complex
with $7$ vertices  indexed by the summands of $R_1$: three vertices
span a $2$-simplex if and only if the orthogonal complement of the sum of their
associated $A_1$-summands is of type $D_4$. The preceding lemma tells us that every edge is in exactly one
$2$-simplex. Probably the $W(E_7)$-stabilizer of $R_1$ is the full
automorphism group of this complex.
\\

The  subsystems of type $7A_1$ make up two $\Scal_7$-orbits, represented by
\begin{enumerate}
\item[(A)] $\la h_7,h_{12},h_{34},h_{56}, h_{127}, h_{347},h_{567}\ra$, $105$ in number and
\item[(B)] $\la h_{123}, h_{145},h_{167}, h_{256}, h_{247}, h_{357}, h_{346}\ra$, of which there are $30$.
\end{enumerate}
We designate by the same letters $(A)$ and $(B)$ the type  of the corresponding product of roots.

\begin{lemma}\label{lemma:s7action}
Let $F$ be a product of roots of type $(B)$. Then the $\Scal_7$-stabilizer
of $F$ acts transitively on its factors and has order $7.3.2^3$. The subgroup
that stabilizes a given factor is isomorphic to $\Scal_4$.
\end{lemma}
\begin{proof}
Let $F$ be of type $(B)$ and let $\alpha$ be a factor of $F$. Write $\alpha=h_{abc}$.
Then the other factors are of the form $h_{axy}$, $h_{azw}$, $h_{bxz}$, $h_{byw}$, $h_{cxw}$, $h_{cyz}$, where $x,y,z,w$ are the distinct elements of
$\{1,2,\dots ,7\}-\{ a,b,c\}$. So these factors are given by an indexing by $a,b,c$ of the three ways we can split $\{1,2,\dots ,7\}-\{ a,b,c\}$ into two pairs. This description
proves that $\Scal_7$ is transitive of the collection of pairs $(F, \alpha)$
with stabilizer mapping isomorphically onto the permutation group of $\{x,y,z,w\}$.
The lemma follows.
\end{proof}

\begin{lemma}\label{lemma:s3property}
The space $\coble_2$ is spanned by the $30$ root products
of type $(B)$ and is annihilated  by $\sum_{w\in \Scal(i,j,k)} \sign(w) w$ for any
$3$-element subset $\{i,j,k\}$ of $\{1,\dots ,7\}$.
\end{lemma}

For the proof we need:

\begin{lemma}\label{lemma:d4}
A root system $S$ of type $D_4$ contains exactly $3$ subsystems
of type $4A_1$  and  the discriminants of these three subsystems
(relative some choice of positive roots) are such that a signed sum is zero.
More precisely, if $S_o\subset S$ is a subsystem of type $4A_1$,
then the $(8)$ reflections in $W(S)-W(S_o)$ decompose into two equivalence classes
with the property that two reflections $s,s'$  belong to different classes if and only if they do not commute. In that case  $S_o)$, $sS_o$, $s'S_o$
are the distinct $4A_1$-subsystems of $S$ and if $f$ is the product of $4$
pairwise perpendicular roots in $S_o$, then $f=s(f)+s'(f)$.
The plane in the fourth symmetric power of the complex span of the root system
generated by these discriminants  affords an irreducible
representation of the Weyl group of the root system.
\end{lemma}
\begin{proof}
In terms of the standard model for the $D_4$-system,  the set of vectors $\pm \eps_i\pm \eps_j$,
$1\le i<j\le 4$ in Euclidean $4$-space, the $4A_1$-subsystems correspond to the three ways of partitioning $\{1,2,3,4\} $ into  parts of size $2$. For instance, the partition
$\{\{1,2\},\{3,4\}\}$ yields  $\{ \pm \eps_1\pm \eps_2,\pm \eps_3\pm \eps_4\}$, whose discriminant is
(up to sign) equal to  $\eps_1^2\eps_3^2+\eps_2^2\eps_4^2-\eps_1^2\eps_4^2+\eps_2^2\eps_3^2$.
We can verify the lemma for $s=s_{\eps_1-\eps_3}$, $s'=s_{\eps_1-\eps_4}$ and deduce the
general case from that.

The last clause is easily verified.
\end{proof}

\begin{proof}[Proof of Lemma \ref{lemma:s3property}]
Consider $F= h_7h_{12}h_{34}h_{56} h_{127} h_{347}h_{567}$ (a typical root product of type $(A)$). The four factors that are not of type $(1)$, $h_7,h_{12},h_{34},h_{56}$,
lie in a subsystem of type $D_4$. If we  let $s$ resp.\ $s'$ be the reflection in
$h_{135}$ resp.\ $h_{246}$, then
\begin{gather*}
s(h_7h_{12}h_{34}h_{56})=h_{246}h_{235}h_{145}h_{136}\\
s'(h_7h_{12}h_{34}h_{56})=h_{135}(-h_{146})(-h_{236})(-h_{245})
=-h_{135}h_{146}h_{236}h_{245}.
\end{gather*}
Notice that the second product is obtained from the first by applying minus the transposition
$(34)$. According to Lemma \ref{lemma:d4} we then have
\begin{equation*}\label{AB}
h_7h_{12}h_{34}h_{56}=
\left(1-(34)\right)h_{246}h_{235}h_{145}h_{136}.\tag{AB}
\end{equation*}
After multiplying both sides with $h_{127}h_{347}h_{567}$, we see
that $F$ has been written as a difference of two products of type $(B)$:
$f=\left(1-(34)\right)G$ with
\[
G:=h_{127}h_{347}h_{567}h_{246}h_{235}h_{145}h_{136}.
\]
In particular, the type $(B)$-products generate $\coble_2$.
It follows from Lemma \ref{lemma:s7action} that the  $\Scal_7$-stabilizer of $G$ has two orbits in the collection of $3$-element subsets $\{i,j,k\}\subset\{1,\dots ,7\}$: those for which $h_{ijk}$ is a factor of $G$  and those for which there exist a factor $h_{abc}$ of $G$
with $\{a,b,c\}\cap\{i,j,k\}=\emptyset$. So there are only
two cases to verify.

We first do the case $I=\{3,4,5\}$. We are then in the second case
because $I\cap \{1,2,7\}=\emptyset$ and $h_{127}$ is a factor of $G$.
To this purpose we look at the $D_4$-system defined by the pair $h_7$,
$h_{127}$:
it is the system that contains the four roots $h_{34}$, $h_{56}$, $h_{347}$, $h_{567}$.
The reflections $s$ resp.\ $s'$ perpendicular to $h_{367}$ resp.\ $h_{467}$ lie in this $D_4$ summand and do not commute.We have
\begin{align*}
s(h_{34}h_{56}h_{347}h_{567})=&h_{467}(-h_{357})h_{64}h_{35}\\
s'(h_{34}h_{56}h_{347}h_{567})=&(-h_{367})(-h_{457})h_{63}h_{45}
\end{align*}
so that
\[
h_{34}h_{56}h_{347}h_{567}-h_{467}h_{357}h_{46}h_{35}+h_{367}h_{457}h_{36}h_{45}=0.
\]
The second resp.\ third term are obtained from the first by
applying to it minus the transposition $(45)$ resp. minus the transposition $(35)$, so that
\[
\left(1-(45)-(35)\right)h_7h_{12}h_{127}h_{34}h_{56}h_{347}h_{567}=0.
\]
If we combine this with \eqref{AB}, and observe that
\[
\left(1-(45)-(35)\right)\left(1-(34)\right)=\sum_{w\in\Scal (3,4,5)} \sign (w)w,
\]
then we find that the latter kills
$G=h_{127}h_{347}h_{567}h_{246}h_{235}h_{145}h_{136}$.

An instance of the first case, namely the assertion that $G$ is also killed by
$\sum_{w\in\Scal (1,2,7)} \sign (w)wG$, follows by exploiting the symmetry
properties of $G$: the transpositions $(34)$ and $(35)$ have the same effect
on $G$ as resp.\ $(12)$ and $(17)$. This implies that
\[
\sum_{w\in\Scal (1,2,7)} \sign (w)wG=\sum_{w\in\Scal (3,4,5)} \sign (w)wG=0.
\]
\end{proof}

\begin{corollary}\label{cor:coble7}
We have $\dim \coble_2\le 15$.
\end{corollary}
\begin{proof}
Let $f\in\coble_2$. Since a monomial of type $(B)$ has a unique factor of the form $(12 a)$, $a\in\{ 3,\dots,7\}$, we can write $f$ accordingly: $f=(124)f_3+\cdots +(127)f_7$. For every $a$, we have $6$ type $(B)$ monomials
corresponding to the ways we index the splittings of
the complement of $a$ in $\{3,4,5,6,7\}$ by $\{1,2,a\}$. The symmetric group
$\Scal (1,2,a)$ permutes these six root products simply transitively. These root products
satisfy the corresponding alternating sum relation  and so we can arrange that each $f_a$ is a linear combination of $6$ monomials whose alternating sum of coefficients is zero.  If we take as our guiding idea to make $a$ as small as possible, then
it turns out that in half of the cases we can do better.

Let us first assume $a\in\{5,6,7\}$. We then invoke the relation  defined by $\{3,4,a\}$:
\[
\sum_{w\in\Scal(3,4,a)}\sign(w)w f_a=0.
\]
Four of the six terms have a factor $(123)$ or $(124)$, whereas the other two
have a factor $(12a)$ and combine to $(1-(34))f_a$. So this relation allows us to
arrange that $f_a$ and $(34)f_a$ have the same
coefficient. We thus make $f_a$ vary in a space of dimension $\le 3$.
If $a\in\{6,7\}$, then we can repeat this game with $\{4,5,a\}$. This allows us to assume in addition that $f_a$ and $(45)f_a$ have the same coefficient. But then $f_a$ must have all its coefficients equal. So $f_a$ varies in a space of dimension $\le 1$ for
$a=6,7$, of dimension $\le 3$ for  $a=5$, and of dimension $\le 5$ for $a=3,4$. This proves that $\dim\coble_2\le 15$.
\end{proof}

\begin{proof}[Proof of Proposition \ref{prop:coble7}]
If we combine Proposition \ref{prop:macdonald} and Corollary \ref{cor:coble7},
we see that $\coble_2$ is an  irreducible representation of $W(R)$ of dimension
$\le 15$. Since $W(R)=W(R)_+\times\{1,c\}$, it  will then also be an irreducible representation of $W(R)_+$. The symmetric bilinear form on the root lattice
induces a nondegenerate form on the root lattice modulo two times the weight lattice
(this is a $\FF_2$-vector space of dimension $6$). This identifies $W(R)_+$
with the symplectic group
$\Sp (6,\ZZ/2)$ and it is well-known (see for instance  \cite{atlas}, where this group
is denoted $S_6(2)$) that the irreducible representations  of dimension $<15$
are the trivial representation, the sign representation and the standard representation of degree $7$. It is easy to see that $\coble_2$ is neither of these.
Since $c$ acts as $-1$ in $\coble_2$, the proposition follows.
\end{proof}

The roots orthogonal to an $A_5$-subsystem of an $E_6$-system make up a
system of type $A_1$ or  $A_2$.  In terms of our root basis, they are
represented by $\la h_{23},h_{34},h_{45},h_{56},h_{67}\ra$ (with $\la h_1\ra$ as
the  perpendicular system) and  the system $\la h_{123},h_{34},h_{45},h_{56},h_{67}\ra$ (with $\la h_1,h_{12}\ra$ as perpendicular system).  We shall call an $A_5$-subsystem of the second type \emph{special}. Conversely,  the roots perpendicular to a $A_2$-subsystem form a special $A_5$-system. Since the $A_2$-subsystems make
up a single Weyl group equivalence class, the same is true for the special $A_5$-subsystems.

\begin{proposition}\label{prop:fixpart}
If we regard $\coble_2$ as a vector space of degree $7$ polynomials on $\hfrak_7$,
then their common zero set  $Z_7$ is the union of the linear subspaces perpendicular
to a system of type $D_4$ or to a special system of type $A_5$.
\end{proposition}
\begin{proof}
The $D_4$-subsystems constitute  a single Weyl group equivalence class and so
we may take as our system the one spanned by the fundamental roots
$\la h_{23}, h_{34}, h_{45}, h_{123}\ra$. We must show that every $7A_1$-subsystem of $R$ meets this $D_4$-system. The positive roots of the $D_4$-system are
$\{ h_{ij}\}_{2\le i<j\le 5}$ and $\{ h_{1ij}\}_{2\le i<j\le 5}$. It is easy to see from our description that every $7A_1$-system of type $(A)$ contains a root $h_{ij}$ with
$2\le i<j\le 5$. Similarly, we see that every $7A_1$-system of type $(B)$ contains a root $h_{1ij}$ with $2\le i<j\le 5$.

We argue for the special $A_5$-system
$\la h_{123},h_{34}, h_{45}, h_{56}, h_{67}\ra$
 in a similar fashion.  Its positive roots are
$\{ h_{ij}\}_{3\le i<j\le 7}$ and $\{ h_{12i}\}_{3\le i\le 7}$.
Every $7A_1$-system of type $(A)$ contains a root $h_{ij}$ with  $3\le i<j\le 7$ and every $7A_1$-system of type $(B)$ contains a root $h_{12i}$ with $3\le i\le 7$.

It remains to show that this exhausts $Z_7$. Every subsystem
of $R$  that does not contain a $D_4$-subsystem has only components of type $A$.
If in addition it does not contain a  special $A_5$-system, then  any such a subsystem is Weyl group-equivalent to a proper subsystem in the (nonsaturated) $A_7$-system spanned  by the fundamental roots
$h_{12},h_{34},\cdots ,h_{67}$ and the highest root $h_1$. The latter has as its positive roots $\{ h_{ij}\}_{1\le i<j\le 7}\cup \{h_i\}_{i=1}^7$ and  is therefore disjoint with
the $7A_1$-subsystem $\la h_{123}, h_{145},h_{167}, h_{256}, h_{247}, h_{357}, h_{346}\ra$. This implies that $Z_7$ is as asserted.
\end{proof}

\begin{question}
Is $\coble_2$ the space of degree $7$ polynomials on $\hfrak_7$ that vanish on $Z_7$?  We expect this to be equivalent to the question whether $\coble_2$ is complete as a linear system on $\Mcal^*_{m,2}$.
\end{question}

\section{The Coble linear system}\label{sect:system}
Let $A$ be a vector space of dimension three so that $\PP(A)$ is a projective plane.
Given  a numbered set $(p_1,\dots,p_N)$ of $N\ge 5$ points in $\PP(A)$ that are in generic position, then for any $5$-tuple $(i_0,\dots ,i_4)$ with $i_0,\dots ,i_4$ pairwise distinct and taken from $\{ 1,\dots ,N\}$, the four ordered  lines
$p_{i_0}p_{i_1}$, $p_{i_0}p_{i_2}$,
$p_{i_0}p_{i_3}$, $p_{i_0}p_{i_4}$ through $p_{i_0}$
have a cross ratio. The collection of cross ratio's thus obtained make up a complete projective invariant of $(p_1,\dots,p_N)$: we may choose coordinates such that $p_1=[1:0:0]$,  $p_2=[0:1:0]$,  $p_3=[0:0:1]$,  $p_4=[1:1:1]$ and the coordinates
$[z_0:z_1:z_2]$ for $p_i$, $i>4$, are then given by cross ratio's.
For instance,  $z_1:z_2=(p_1p_2:p_1p_3:p_1p_4:p_1p_i)$. If $a_i\in A$
represents $p_i$, then we can write the this as a cross ratio
of $4$ lines in the plane $a_1\wedge A\subset \wedge^2A$:
$(a_1\wedge a_2:a_1\wedge a_3:a_1\wedge a_4:a_1\wedge a_i)$.

Now let us observe that if $(v_1,v_2,v_3,v_4)$
is a generic ordered $4$-tuple in a vector space $T$ of dimension two, then
the corresponding points  in $\PP(T)$ have a cross ratio that can be written as a ratio of two elements of $\det (T)^2$, namely
$(v_1\wedge v_4)(v_2\wedge v_3): (v_2\wedge v_4)(v_1\wedge v_3)$.
If we apply this to to the present case, then we get
\[
z_1:z_2=(p_1p_2:p_1p_3:p_1p_4:p_1p_i)=|12i| |134|: |13i| |124|,
\]
where we used the Coble notation. Thus the cross ratio's formed in this manner
allow us to reconstruct $(p_1,\dots ,p_5)$ up to projective equivalence.
We can express this in terms of roots as follows.

\begin{lemma}\label{lemma:rootcrossratio}
Let $(S; e_1,\dots ,e_{9-d})$ be a marked Del Pezzo surface of degree $d\le 4$,
$S\to \check{\PP}^2(H^0(S,\ell))$ the contraction morphism defined by the linear system
$|\ell|$ (as usual) and $p_i$ the image of $E_i$. Then
$p_1,\dots ,p_4$ are in general position if and only if none of the roots
in the $A_4$-subsystem generated by $(h_{123},
h_{12},h_{23},h_{34})$
is effective in $\pic (S)$.

If that is the case and  $K$ is a cuspidal anticanonical curve on $S$, then
for  $i>4$ the cross ratio $(p_1p_2:p_1p_3:p_1p_4:p_1p_i)$ equals the ratio
of the two elements of the line $\pic (K)^o$ given by
$r_K(h_{2i}h_{12i} h_{34}h_{134})$ and $r_K(h_{3i}h_{13i} h_{24}h_{124})$.
\end{lemma}
\begin{proof}
The first part is left to the reader as an exercise.
As to the second part, choose affine coordinates $(x,y)$ in $\check{\PP}^2(H^0(S,\ell))$ such that the $K$ is given
by $y^3=x^2$. So $p_i=(t_i,t_i^3)$
for some $t_i$. For $i\not=1$, the line $p_1p_i$ has tangent
$[t_i-t_1:t_i^3-t_1^3]=[1: t_i^2+t_it_1+t_1^2]$. So the cross ratio  of the lines $p_1p_i$
involves  factors of the form
\[
(t_j^2+t_jt_1+t_1^2)-(t_i^2+t_it_1+t_1^2)=(t_j-t_i)(t_j+t_i+t_1), \quad 2\le i<j\le 5.
\]
If we use the $x$-coordinate to identify $\pic_0(K)$ with $\CC$,
then such a factor can be written $r_K(h_{ij}h_{1ij})$. The last assertion follows.
\end{proof}

\begin{remark}
Notice that the roots that appear in the numerator resp.\ the denominator of
$h_{25}h_{12i} h_{34}h_{134}: h_{3i}h_{13i} h_{24}h_{124}$  are four pairwise
perpendicular roots which all lie in a single $D_4$-subsystem.
\end{remark}

\subsection*{The Coble system in the  degree four case}
We first consider  a Fano surface of degree $5$. We recall that such a surface
$S$ can be obtained by blowing  up $4$ points of a projective plane in general position,  and so is unique up to isomorphism. Any automorphism of this surface that acts trivially on its Picard group preserves every exceptional curve ($=$ line) and hence is the identity. It follows that the automorphism group of $S$ is the Weyl group $W(A_4)$. There  are 10 lines on $S$. If five of them
make up pentagon, then their sum is an anticanonical divisor. There are 12 such pentagons and they generate the anticanonical system. Let us now fix a marking
$(e_1,e_2,e_3,e_4)$ for $S$. To every $p\in S$ we associate a marked Del Pezzo surface $(S_p;e_1,\dots ,e_5)$ of degree $4$ by letting $S_p$ be the  blowup of $S$ in $p$ and letting $e_5$ be the class of the exceptional divisor. This defines a rational map $S\dashrightarrow \Mcal^*_{m,4}$.  We will see that this is in fact an isomorphism.

\begin{proposition}\label{prop:coble5}
The $12$ Coble covariants for $D_5$ span a complete linear system of dimension
$6$ and define an embedding of  $\Mcal_{m,4}^*$ in a projective space of dimension
$5$. The image is `the' anticanonically embedded Fano surface  of degree $5$
(so that the Coble system is anticanonical) with $\Mcal_{m,4}^*-\Mcal_{m,4}^\circ$
mapping onto the union of its ten lines. The divisor of every Coble covariant
is a pentagon on this Fano surface and every pentagon thus occurs.
\end{proposition}
\begin{proof}
Let $A$ be a complex vector space of dimension $3$ and let  $p_1,\dots ,p_5\in \PP(A)$.
We first assume that $p_1,\dots ,p_4$ are in general position (i.e., no three collinear). We then adapt our coordinate system accordingly: $p_i=[a_i]$, with
$a_1=(1,0,0)$, $a_2=(0,1,0)$, $a_3=(0,0,1)$ and $a_4=(1,1,1)$. If $a_5=(z_0,z_1,z_2)$, then  typical determinants involving $a_5$ are:
\[
|1 2 5|=z_2, \quad |14 5|=z_2-z_1.
\]
So for instance
\[
 |415| |152| |523| |234| |341|=(z_1-z_2).(-z_2).z_0.1.-1=z_0z_1z_2-z_1^2z_2.
 \]
We obtain in this  manner the $6$  polynomials $z_0z_1z_2-z_i^2z_j$, $i\not=j$,
and it is easily verified that any other  Coble covariant is a linear combination of these.
They are visibly linearly independent and hence form a basis for $\coble_5$. It is
precisely the linear system of cubic curves that
pass through $p_1,\dots ,p_4$. So the Coble system is anticanonical  and defines an embedding of the blowup $S$ of $\PP(A)$ in $p_1,\dots ,p_4$ to $\PP^5$.
The remaining assertions are verified in a straightforward manner.
\end{proof}

\begin{remark}
This proposition and its proof show that the moduli space $\Mcal^*_{m,4}$
is as a variety simply obtained from $\PP^2$ by blowing up the vertices
$p_1,\dots ,p_4$  of the coordinate simplex. This argument then also shows
that there is universal semistable marked Del Pezzo surface of degree $4$,
$\Scal_{m,4}\to\Mcal^*_{m,4}$:

over $p_5\in \Mcal^*_{m,4}$ we put the blowup of the surface $\Mcal^*_{m,4}$ in $p_5$ so that  $\Scal_{m,4}$ is simply $\Mcal^*_{m,4}\times\Mcal^*_{m,4}$ blown up along the
diagonal with one of the projections serving as the structural morphism.
\end{remark}

\subsection*{The Coble system in the  degree three case}
Our discussion starts off with the following lemma.

\begin{lemma}\label{lemma:nobasept3}
The  Coble system $\coble_3$ has no base points.
\end{lemma}
\begin{proof}
The Coble linear system pulled back to $\PP(\hfrak_6)$ has according to
Proposition \ref{prop:z6} as its base point locus
the projective arrangement $\PP(Z_6)$, the union of the fixed point  hyperplanes
of Weyl subgroups of type $A_3$. Since $\Mcal_3^*$ is a quotient of
$\PP(\hfrak_6)-\PP(Z_6)$ it follows that $\coble_3$ has no base points
and hence defines a morphism to a $\PP^9$.
\end{proof}

We use what we shall call the \emph{Naruki model} of $\Mcal^*_{m,3}$.
This is based on a particular way of getting a degree $3$ Fano surface as a
blown-up projective plane: we suppose the points in question to be labelled
$p_i,q_i$ with $i\in\ZZ/3$ and to lie on the coordinate lines of $\PP^2$ as follows:
\[
\begin{matrix}
p_0=[0:1:a_0], &  p_1=[a_1:0:1],& p_2=[1:a_2:0],\\
q_0=[0:1:b_0], &  q_1=[b_1:0:1],& q_2=[1:b_2:0].
\end{matrix}
\]
For the moment we assume that blowing up these points gives rise
to a Fano surface $S$ so that in particular none of the $a_i, b_i$ is zero
and $a_i\not=b_i$. If we blow up these points,
the strict transform of the coordinate triangle is  an anticanonical curve  $K$ (it is a tritangent of the corresponding cubic surface). This `partial rigidification' reduces the projective linear group
$\PGL(3,\CC)$ to its maximal  subtorus that leaves the  coordinate triangle invariant.
The following expressions are invariant under that torus
\[
\alpha_i:=a_i/b_i\, ( i\in \ZZ/3),\quad  \delta:=-b_0b_1b_2,
\]
and together they form a complete projective invariant of the configuration. Notice the formulas
\[
a_0b_1b_2=-\alpha_0\delta, \quad a_0a_1b_2=-\alpha_0\alpha_1\delta,\quad
a_0a_1a_2=-\alpha_0\alpha_1\alpha_2\delta.
\]
As explained in the Appendix of \cite{naruki}, $\alpha_i$ and $\delta$ have a simple interpretation in terms of the of pair $(S,K)$: if we denote the exceptional curves $A_i ,B_i$, and
$E$ is the strict transform of the line through $q_1$ and $q_2$, then the cycles
$A_i-B_i$, $ i\in \ZZ/3$ and $B_0-E$ span in $\pic (S)$ the orthogonal complement to the components of $K$ and the numbers in question can be interpreted as
their images in $\pic^0(K)\cong\CC^\times$. The classes themselves make up the basis of a $D_4$ root system with the last one representing the central node. We denote the $4$-torus for which $\alpha_0,\alpha_1,\alpha_2 ,\delta$ is a basis of characters by $T$. This torus comes with an action of $W(D_4)$ and this makes it an adjoint torus of type $D_4$. We denote by $T^\circ$ the open set of its regular elements. This is the complement of the union of reflection hypertori, i.e., the locus where none of the $D_4$-roots is $1$. It has the interpretation as the moduli space of marked nonsingular cubic surfaces with the property that a particular tritangent (which is entirely given by the marking) has not its three lines collinear. If we want to include that case too, we must first blow up the identity element of $T$, $\bl_1(T)\to T$, and then remove the strict transforms of the reflection hypertori. This open subset, $\bl_1(T)^\circ\subset \bl_1(T)$, is a model for the moduli space of marked nonsingular cubic surfaces, in other words, it can be identified with $\Mcal^\circ_{m,3}$.  The modular interpretation implies that this variety has  a  $W(E_6)$-action, although  only the action of a Weyl subgroup of type $D_4$ is manifest. The action of the two missing fundamental reflections was written down by Naruki and Sekiguchi: the one that in the Dynkin diagram is attached to $\alpha_i$ is given
by
\begin{align*}
\alpha_i&\mapsto -\alpha_0\alpha_1\alpha_2\delta^2\frac{1-\alpha_i}{
1-\alpha_0\alpha_1\alpha_2\delta^2}\\
\alpha_{i\pm 1}&\mapsto\frac{(1-\alpha_0\alpha_{i\pm 1}\delta)(1-\alpha_0\alpha_1\alpha_2\delta)}{(1-\alpha_0\delta)(\alpha_{i\pm 1}-\alpha_0\alpha_1\alpha_2\delta)}\\
\delta&\mapsto\delta^{-1}\frac{(1-\alpha_i\delta)(1-\alpha_0\alpha_1\alpha_2\delta^2)}{( 1-\alpha_i)(1-\alpha_0\alpha_1\alpha_2\delta)}.
\end{align*}
(We give this formula because of its remarkable form only---we shall not use it.)
Subsequently Naruki \cite{naruki} found a nice $W(E_6)$-equivariant
smooth
projective completion of this space with a normal crossing divisor as boundary. What
is more relevant  here is a (projective) blow-down of  his completion that was also introduced by him. We shall take Naruki's construction of the  latter as our guide and reprove some of his results in the process.

We identify the complexification of $\Hom (\CC^\times,T)$  with the Lie algebra $\tfrak$ of $T$, so that the latter has a natural $\QQ$-structure.   The decomposition $\Sigma$ of
$\tfrak(\RR)$ into its $W(D_4)$-chambers has its rays spanned by the coweights
that  lie in the $W(D_4)$-orbit of a fundamental coweight. One of these is the orbit
of coroots and has  24 elements; the other three consist of minuscule weights and are a single orbit under the full automorphism group of the $D_4$-system;
it has also 24 elements. If we remove the faces that contain a coroot we obtain a coarser decomposition of $\tfrak(\RR)$ that we denote by $\Sigma'$;
a maximal  face of $\Sigma'$ is now an orbit of a Weyl chamber under the stabilizer
of a coroot (a type $3A_1$-Weyl group).

Let $T\subset T_\Sigma$ be the associated torus embedding. It is a smooth with normal crossing boundary. The boundary divisors are in bijective correspondence with the above coweights. We shall refer to those that correspond to coroots resp.\ minuscule weights as \emph{toric coroot divisors} resp.\ \emph{toric minuscule weight divisors}. So there are 24 of each.

Now blow up successively in $T_\Sigma$: the identity element (in other words the fixed point set of $W(D_4)$), the fixed point sets of the Weyl subgroups of type $A_3$, the fixed point sets of the Weyl subgroups of type $A_2$. We denote the resulting blowup $\hat T_\Sigma$. In this blowup  the exceptional divisors of type $A_3$ have been separated and each is  naturally a product. To be precise,  a $W(A_3)$-Weyl subgroup  $G\subset W(D_4)$ has as its fixed point locus in $\hat T_\Sigma$ a copy of a $\PP^1$ whose tangent line at the identity is the $G$-fixed point set in $\tfrak$ and the divisor associated to $G$ is then naturally the product of the projective line  $\hat T^{G}_\Sigma$ and the projective plane $\PP (\tfrak/\tfrak^G)$ blown up in the fixed points of the $A_2$-Weyl subgroups of $G$ (there are four such and they are in general position).

The Coble system on Naruki's completion, together with the $W(E_6)$-action on it, was identified in \cite{geemen} (5.9 and 4.5). It
is the pull-back of
$$
H^0(T_\Sigma,\Ocal(D_S+2D_R)\otimes \mfrak_e^3)
$$
to $\hat T_\Sigma$, here $D_S$, $D_R$ are the sum of the 24 divisors
corresponding to the rays spanned by the minuscule weight and the coroots
respectively, and $\mfrak_e$ is the ideal sheaf of the identify element
$e\in T_\Sigma$.

We now give an explicit description of
the Coble covariants in terms of $(\alpha_0,\alpha_1,\alpha_2,\delta)$.
For this we begin with observing the following simple identities:
\[
|p_0p_1p_2|= a_0a_1a_2+1,\quad |p_iq_ip_{i+1}|=(b_i-a_i)a_{i+1},\quad |p_iq_ip_{i-1}|=b_i-a_i.
\]
A straightforward computation yields
\[
|p_0p_1p_2q_0q_1q_2|=\pm (b_0-a_0)(b_1-a_1)(b_2-a_2)(1-a_0a_1a_2b_0b_1b_2).
\]
We substitute these values in the formulae for the Coble covariants, but for reasons
similar as in Section \ref{sect:cusp}  we divide these by $(b_0-a_0)(b_1-a_1)(b_2-a_2)$.  For example,
\begin{multline*}
|p_0p_1p_2q_0q_1q_2|.|p_0q_0p_1|.|q_1p_2q_2|=\\
=\pm (1-a_0a_1a_2b_0b_1b_2).(b_0-a_0)a_1.(b_2-a_2)\\
=\pm \alpha_1\delta (1-\alpha_0\alpha_1\alpha_2\delta)(1-\alpha_0)(1-\alpha_2)
\end{multline*}
and
\begin{multline*}
|p_0q_0p_1|.|p_0q_0q_2|.|p_1q_1q_2|.|p_1p_2q_2|.|p_0q_1p_2|.|q_0q_1p_2|=\\
=\pm \frac{(b_0-a_0)a_1. (b_0-a_0). (b_1-a_1)b_2.(b_2-a_2)(a_0b_1a_2+1)
(b_0b_1a_2+1)}{(b_0-a_0)(b_1-a_1)(b_2-a_2)}\\
=\pm\alpha_0\alpha_1\delta (1-\alpha_0)(1-\alpha_0\alpha_2\delta)(1-\alpha_2\delta).
\end{multline*}
We thus find  for the Coble covariants (40 up to sign) the following expressions:
\[
\begin{tabular}{cc}
$1$&$\delta(1-\alpha_0)(1-\alpha_1)(1-\alpha_2)$\\
$2$&$\alpha_0\alpha_1\alpha_2\delta^2(1-\alpha_0)(1-\alpha_1)(1-\alpha_2)$\\
$3$&$(1-\alpha_0\delta)(1-\alpha_1\delta)(1-\alpha_2\delta)$\\
$4$&$\alpha_0\alpha_1\alpha_2\delta(1-\alpha_0\delta)(1-\alpha_1\delta)(1-\alpha_2\delta)$\\
$5$&$(1-\alpha_0\alpha_1\delta)(1-\alpha_1\alpha_2\delta)(1-\alpha_2\alpha_0\delta)$\\
$6$&$\delta(1-\alpha_0\alpha_1\delta)(1-\alpha_1\alpha_2\delta)(1-\alpha_2\alpha_0\delta)$\\
$7_i$&$(1-\alpha_{i-1}\delta)(1-\alpha_{i+1}\delta)(1-\alpha_0\alpha_1\alpha_2\delta)$\\
$8_i$&$\alpha_i\delta(1-\alpha_{i-1})(1-\alpha_{i+1})(1-\alpha_0\alpha_1\alpha_2\delta^2)$\\
$9_i$&$\delta(1-\alpha_{i-1})(1-\alpha_{i+1})(1-\alpha_0\alpha_1\alpha_2\delta^2)$\\
$10_i$&$\alpha_{i-1}\alpha_{i+1}\delta (1-\delta)(1-\alpha_i\alpha_{i-1}\delta)(1-\alpha_i\alpha_{i+1}\delta)
$\\
$11_i$&$(1-\delta)(1-\alpha_i\alpha_{i+1}\delta )(1-\alpha_i\alpha_{i-1}\delta)$\\
$12_i$&$\alpha_i\delta(1-\alpha_{i+1}\delta)(1-\alpha_{i-1}\delta)(1-\alpha_0\alpha_1\alpha_2\delta)$
\\
\\
$13$&$(1-\delta)(1-\alpha_0\alpha_1\alpha_2\delta)(1-\alpha_0\alpha_1\alpha_2\delta^2)$\\
$14_i$&$\alpha_{i-1}\alpha_{i+1}\delta (1-\delta)(1-\alpha_i)(1-\alpha_i\delta)$\\
$15_i$&$(1-\alpha_i\delta)(1-\alpha_{i-1}\alpha_{i+1}\delta)(1-\alpha_0\alpha_1\alpha_2\delta^2)$\\
$16_i$&$\delta(1-\alpha_i)(1-\alpha_{i-1}\alpha_{i+1}\delta)(1-\alpha_0\alpha_1\alpha_2\delta)$\\
$17_{ijk}$&$\alpha_i\delta (1-\alpha_j)(1-\alpha_k\delta)(1-\alpha_j\alpha_k\delta)$\\
\end{tabular}
\]
Notice that the zero sets  in $T$ of any of these expressions is the union of the reflection hypertori of a 
subsystem of type $3A_1$ (in the first $12$ cases) or the $A_2$ (for the last $5$).
We now quote from  \cite{naruki} (Proposition 11.3):

\begin{theorem}[Naruki]\label{thm:crossratiovar}
There is a projective contraction
\[
\hat T_\Sigma\to \check T_\Sigma
\]
which contracts each $A_3$-divisor along the projection on its $2$-dimensional factor.
The contracted variety $\check T_\Sigma$ is nonsingular and
the action of $W(E_6)$ on $(\bl_e T)^\circ$ extends regularly to it. This action is transitive on the collection of 40 divisors that are of toric coroot type or of $A_2$-type.
\end{theorem}

The 40 divisors in question are easily seen to be pairwise disjoint. Naruki also shows
(Section 12 of \emph{op.\ cit.}) that each of these can be contracted to a point. We
can see that quickly using the theory of torus embeddings:
the 24 toric coroot divisors get contracted if we replace in the above discussion the decomposition $\Sigma$ of  $\tfrak(\RR)$ by the coarser one, $\Sigma'$,  that we obtain by removing
the faces that contain a coroot. The $W(E_6)$-action and the theorem above imply that this contraction is  then also possible for the remaining  16 divisors. The singularities thus created, for instance
the one defined by the ray spanned by the coroot $\delta^\vee$, can be understood as follows:
the  natural affine $T$-invariant neighborhood  of the point of $T_{\Sigma'}$ defined by the ray spanned by the coroot $\delta^\vee$ is $\spec$ of the algebra generated by the elements of the orbit of $\delta$ under the Weyl group of the $3A_1$-subsystem $\la \alpha_0,\alpha_1,\alpha_2\ra$: $\spec\CC[\delta,\alpha_0\delta,,\alpha_1\delta, \dots ,\alpha_0\alpha_1\alpha_2\delta]$.
This is a cone over the Veronese embedding  of $(\PP^1)^3$.

The following more precise result is in \cite{geemen} (Theorem 5.7) and
follows also from \cite{yoshida}.

\begin{theorem}
The Coble covariants generate on $\check T_\Sigma$ a linear system without base points
that has the property that its restriction each of the 40 divisors of  toric coroot type or of $A_2$-type is trivial. The resulting morphism to a nine dimensional projective space realizes Naruki's contraction.
\end{theorem}
\begin{proof}[Outline of proof]
We have a natural decomposition of $T_{\Sigma}$  into strata by type: $D_4$ (yielding the identity element), $A_3$, $A_2$ and $\{ 1\}$ (being open in
$T_\Sigma$). There is a corresponding decomposition of $\hat T_\Sigma$ (and of
$\check T_\Sigma$, but we find it more convenient to work on the former), albeit that strata are then indexed by chains of strata in $T_\Sigma$ that are totally ordered for incidence.
We first check that along every stratum the Coble covariants define, modulo the stated contractions property, an embedding.  For this, the $W(D_4)$-equivariance allows us to concentrate on the open
subset $U=\spec( \CC[\alpha_0,\alpha_1,\alpha_2,\delta])$ of $T_\Sigma$ and its preimage
$\hat U$ in $\hat T_\Sigma$.

It is clear from the expressions we found that the Coble covariants generate
$\CC[\delta,\alpha_0\delta,\dots ,\alpha_0\alpha_1\alpha_2\delta]$ after we localize away from
the kernels of the roots (that is, we make each expression $\text{root}-1$ invertible). So we have an embedding on the corresponding  open subset of $\check T_\Sigma$ (this contains the singular point defined the coroot  $\delta^\vee$). A closer look at the equations shows
that this is in fact even true if we allow some of the roots to be $1$, provided that they are mutually perpendicular. In other words, the linear system defines an embedding on the intersection of the open stratum with $U$.

Now let $Z$ be the  $A_3$-stratum that is open in $\alpha_0=\alpha_1=\delta=1$. All Coble covariants vanish on $Z$ and we readily verify that Coble covariants generate the ideal defining
$Z$. Thus the system has no base points on the blowup of $Z$ and the system contracts this
exceptional divisor along the $Z$-direction.

Now let us  look at an $A_2$-stratum $Z'$, say the one that is open in
$\alpha_0=\delta=1$. We observe that the restriction of every covariant to $Z'$ is proportional
with $(1-\alpha_1)(1-\alpha_2)(1-\alpha_1\alpha_2)$ (and can be nonzero). So the linear system
will define the constant map on $Z'$ (or the exceptional divisor over $Z'$).

We turn to the situation at the identity of $T$. Every covariant vanishes there with order three and has for initial part a product of three roots, viewed as linear forms on $\tfrak$. Up to sign, the roots in such a product are the positive roots of a $3A_1$-subsystem (in the first 12 cases)  or of a $A_2$-system (the last $5$ cases). The preceding  implies that the linear system restricted to a  $A_2$-line in $\PP(\tfrak)$ is constant. With some work  we find that the Coble linear system is without base points and generates the ideal defining the identity away from the union of the $A_2$-loci.

The  remaining strata on $\hat T_\Sigma$ are defined by `flags': chains of strata above totally ordered by incidence, with $\{ 1\}<Z'<Z$ as a typical degenerate case.
In that situation, one checks that the covariants generate on that stratum the ideal
$\Ical_{\{ e\}}\Ical_{Z'}\Ical_Z$. We thus see that we have a local embedding along this stratum. The other strata are dealt with in the same way.

This shows that the linear system defines local embeddings modulo the contraction property. So $\check{T}_\Sigma$ is defined as a projective quotient of $\hat T_\Sigma$
and the linear system maps it as local embedding to a nine dimensional projective space. It remains to see that the images of the strata are disjoint. This is left to the reader.
\end{proof}

\begin{remark}
The construction of the  Naruki quotient comes with a stratification  and as one may expect,
each of its members  has a modular interpretation. We here give that interpretation
without proof. As mentioned, the open stratum $T^\circ$ is the moduli space of systems
$(S;e_1,\dots ,e_6;K)$ with $(S;e_1,\dots ,e_6)$ a marked Fano surface of degree $3$ (equivalently, a cubic surface) and $K$ an anticanonical  divisor made up of three nonconcurrent exceptional curves (so that the isomorphism $S\to \bar S$ maps $K$ onto a tritangent $\bar K\subset \bar S$).  Suppose now that $S$ is merely a Del Pezzo surface whose  configuration of $(-2)$-curves is nonempty, but disjoint with $K$.  Then that configuration is of type $rA_1$ and
we are on a stratum contained in $T$ of type $rA_1$ ($1\le r\le 4$) or it is of type $rA_1+sA_1$  with $r\ge 1$ and we are on one of the 24 points that are images of toric coroot divisors (these are the
punctual strata of $T_{\Sigma'}$). In these cases $\bar K$ is a genuine tritangent of $\bar S$
(that lies in the smooth part of $\bar S$). The  other strata are loci for which
$\bar K$ is no longer a tritangent: if $\bar K$ defines an Eckardt point (so that
$K$ consists of three distinct concurrent exceptional curves), then we find ourselves in the stratum that is open in the preimage of the unit element of $T$. If $\bar K$
becomes a union of a double line and another line, then it contains two distinct $A_1$-singularities of $\bar S$ and $K$ is of the form $(2E+E'+C+C'$, where $E,E'$ are exceptional curves
and $C,C'$ are $(-2)$-curves with $E',C,C'$ pairwise disjoint and meeting $E$ normally. We are
then on a stratum that is open in the image of an $A_3$-locus in $\check T$. If $\bar K$ becomes a triple line, then it contains two distinct $A_2$-singularities of $\bar S$ and $K$ is of the form $(3E+C+C'$, where $E$ is exceptional curve and  and $C, C'$ are disjoint $A_2$-curves meeting $E$ normally. We are then representing one of the 16 punctual strata  that are images of an $A_2$-locus in $T$.
\end{remark}

\begin{corollary}
The GIT completion of the moduli space of marked cubic surfaces, $\Mcal^*_{m,3}$,
is $W(E_6)$-equivariantly isomorphic to the Naruki  contraction of $\hat T_\Sigma$.
The Coble linear system embeds $\Mcal^*_{m,3}$ in projective nine space.
\end{corollary}
\begin{proof}
The Coble linear system  on $\hat T_\Sigma$  is without base points and
so the  resulting morphism $f: \hat T_\Sigma\to \PP^9$ realizes the Naruki contraction.
Recall that we have an identification of $\Mcal^\circ_{m,3}$ with $(\bl_e T)^\circ$. This isomorphism  clearly extends to a morphism $\Mcal^*_{m,3}\to f(\hat T_\Sigma)$. This morphism is
birational and since $f( \hat T_\Sigma)$ is normal, it must be a contraction.
\end{proof}

\begin{remark}\label{rem:modularforms}
It is known \cite{act} that  the moduli space of stable cubic surfaces is
Galois covered  by the complex $4$-ball with an arithmetic group $\G$ as Galois group.
The group $\G$ has a single cusp and this cusp represents the minimal strictly stable orbit  of cubic surfaces (i.e., those having three $A_2$-singularities). This gives
$\Mcal_3$ the structure of an arithmetic ball quotient for which $\Mcal^*_3$ is
its Baily-Borel compactification. Allcock and Freitag \cite{af} have used $\G$-modular forms to construct an embedding of this Baily-Borel compactification in a $9$-dimensional projective space. This is precisely the embedding that appears here (see
also Freitag \cite{freitag2} and van Geemen \cite{geemen}).
Via this interpretation it also follows
that the Coble system is complete \cite{freitag}.
\end{remark}

\begin{remark}
The linear system $\coble_3$ can be also studied by restricting it,
as was done in \cite{elbertd5}, to the exceptional divisor
$\PP(\hfrak_5)$ of the blowup of the $e_6$-point in
$\PP(\hfrak_6)$. The generic point of $\PP(\hfrak_5)$ has a
modular interpretation: it parameterizes marked cubic surfaces
with a point where the tangent space meets the surface in the
union of a conic and a line tangent to that conic. The marking
determines the line, but not the conic, for the system of conics
on a cubic surface that lie in plane that contains a given line on that surface 
has two members that are tangent to the  line. So we have
a natural involution  $\iota$ on that space. The projective
space $\PP(\hfrak_5)$ can be seen as the projective span of the
$D_5$-subsystem $R_5$, spanned by the roots not involving $e_6$.

In order to be explicit we also use the standard model for the
$D_5$ root system, i.e., the model for which
$\eps_i-\eps_{i+1}=h_{i,i+1}$ ($i=1,\dots,4$) and
$\eps_4+\eps_5=h_{123}$. This makes $W(R_5)$ the semidirect
product of the group of permutations of the basis elements
$\eps_1,\dots ,\eps_5$ and the group of sign changes in the basis
elements $(\eps_1,\dots ,\eps_5)\mapsto (\pm\eps_1,\dots
,\pm\eps_5)$ with an even number of minus signs. We denote the
basis dual to  $(\eps_1,\dots ,\eps_5)$ by  $(x_1,\dots ,x_5)$.
The Coble covariant we attached to the $3A_2$-system $\la
h_{ij},h_{jk}, h_{lm},h_{m6},h_{ijk}, h\ra$ gives, after dividing
by a common degree 4 factor, the quintic form on $\hfrak_5$ defined by
\[
h_{ij}h_{jk}h_{ik}h_{lm}h_{ijk}=(x_i-x_j)(x_j-x_k)(x_i-x_k)(x_l^2-x_m^2),
\]
whereas the Coble covariant attached to 
$\la h_{ij},h_{ilm},h_{lm},h_{km6},h_{k6}, h_{ijk}\ra$ gives 
\[
h_{ij}h_{jlm}h_{ilm}h_{lm}h_{ijk}=(x_i-x_j)(x_j+x_k)(x_i+x_k)(x_l^2-x_m^2).
\]
These all lie in a single $W(R_5)$-orbit as predicted by Lemma \ref{lemma:d5trans}.
It can be easily checked that the base locus of the system on $\PP (\hfrak_5)$ is the
set of points fixed by a Weyl subgroup of type $A_3$. There are two orbits of subroot systems 
of type $A_3$: one has 40 elements and is represented by 
$\la \eps_1-\eps_2, \eps_2-\eps_3,\eps_3-\eps_4\ra$ and the other has 10 elements and is 
represented by $\la \eps_1-\eps_2, \eps_2-\eps_3,\eps_2+\eps_3\ra$. So the base locus is a 
union of 50 lines. The locus where two such lines meet are the (16) fixed points of a Weyl 
subgroup of type $A_4$ and the (5) fixed points of a Weyl subgroup of type $D_4$.
Blowing up first the 21 points and then the strict
transforms of the 50 lines we obtain a smooth fourfold
$\widetilde{\PP(\hfrak_5)}$ in which the  strict transforms of the planes defined by the (40) 
root subsystems of type $A_2$ have become disjoint. 
The Coble system defines a morphism
\[
\Psi:\,\widetilde{\PP(\hfrak_5)}\longrightarrow \PP^9
\]
which is generically two to one: it identifies the orbits of the involution $\iota$ which 
on $\PP(\hfrak_5)$ is given as the rational map
\[
[x_1:...:x_5]\longmapsto
[x_1^{-1}:...:x_5^{-1}].
\] 
The morphism $\Psi$ is ramified along the exceptional divisors over the $A_4$-points and contracts the
exceptional divisors over the $A_3$-lines to planes and the 40 planes of type $A_2$
to points.
\end{remark}

\subsection*{The Coble system in the  degree two case} An analogue
of Lemma \ref{lemma:nobasept3} holds:

\begin{proposition}\label{lemma:nobasept2}
The linear system $\Ccal_2$ is without base points on $\Mcal^*_{m,2}$. Its restriction
to $\Mcal^\circ_{m,2}$ is an embedding.
\end{proposition}
\begin{proof}
The proof of the first assertion only differs from the one of Lemma \ref{lemma:nobasept2} essentially by replacing the reference to Proposition
\ref{prop:z6} by a reference to Proposition \ref{prop:fixpart}: the pull-back of $\coble_2$
to $\PP(\hfrak_7)$ has according to Proposition \ref{prop:fixpart} as base locus
$\PP(Z_7)$ the projective arrangement defined by the subsystems of type $D_4$
and the special subsystems of type $A_5$. Since the map from
$\PP(\hfrak_7)-\PP(Z_7)$ to $\Mcal_2^*$ is a surjective morphism,  $\coble_2$
is without base points.

The second assertion follows from Lemma \ref{lemma:7a1}:
if we are given a $D_4$-subsystem of a root system of type $E_6$, then
orthogonal to it we have a $3A_1$-system. Thus two disjoint $4A_1$-subsystems
of the given $D_4$-system have a discriminants whose quotient is a quotient of Coble covariants.
According to Lemma \ref{lemma:rootcrossratio} the quotient of two such $4A_1$-subsystems is a cross ratio. Hence all (generalized) cross ratio's are recovered from the Coble covariants. If the seven points are in general position, then these cross ratio's determine the points up to projective equivalence.
\end{proof}

Recall from Proposition \ref{prop:detc2system} that we identified $\coble_2$ with a space of
sections of a square root of $\Ocal_{\Mcal_2^*}(1)$.

\begin{conjecture}\label{prop:c2}
The linear system  $\coble_2$  is without base points and hence defines an injective  morphism
from the moduli space $\Mcal_2^*$ of semistable  quartic curves with level two
structure to a $14$-dimensional  projective space.
\end{conjecture}

We also expect that there is an analogue of the results of Naruki and Yoshida with the role of
$D_4$ taken by $E_6$.


\begin{thebibliography}{99}

\bibitem{af}
D.~Allcock, E.~Freitag:
\textit{Cubic surfaces and Borcherds products},
Comment.\ Math.\ Helv.\ 77  (2002), 270--296.

\bibitem{act}
D.~Allcock, J.A.~Carlson, D.~Toledo:
\textit{The complex hyperbolic geometry of the moduli space of cubic surfaces},
J.\ Algebraic Geom.\  11  (2002),  659--724.

\bibitem{bourbaki:lie} N.~Bourbaki:
\textit{Groupes et alg\`ebres de Lie},
Ch.\ 4,5 et 6, $2^{\text{i\`eme}}$ \'ed., Masson, Paris 1981.



\bibitem{coblebook}
A.B.~Coble:
\textit{Algebraic Geometry and Theta Functions}, AMS New York (1929).

\bibitem{elbert}
E.~Colombo, B.~ van Geemen:
\textit{The Chow group of the moduli space of marked cubic surfaces,}
Ann.\ Mat.\ Pura Appl.\ (4) 183 (2004), 291--316.

\bibitem{elbertd5}
E.~Colombo,  B.~van Geemen:
\textit{A family of marked cubic surfaces and the root system $D_5$,}
available at \texttt{math.AG/0509561}, to appear in Internat. J. Math.


\bibitem{atlas}
J.H.~Conway,  R.T.~Curtis, S.P.~Norton,  R.A.~Parker, R.A.~ Wilson:
\textit{Atlas of finite groups}, Oxford University Press, Eynsham, 1985.


\bibitem{CD}
F.R.~Cossec, I.V.~Dolgachev:
\textit{Enriques Surfaces I}, Birkh\"auser, PM {\bf 76} (1989).

\bibitem{demazure}
M.~Demazure:
\textit{Surfaces de Del Pezzo}, in
S\'eminaire sur les Singularit\'es des Surfaces,
Lecture Note in Math.\ 777, Springer, Berlin-New York, 1980.

\bibitem{dolgort}
I.~Dolgachev, D.~Ortland:
\textit{Point sets in projective spaces
and theta functions},  Ast\'erisque  165 (1988), 210 pp.

\bibitem{dgk}
I.~Dolgachev, B.~van Geemen, S.~Kond\=o:
\textit{A complex ball uniformization of the moduli space of
cubic surfaces via periods of $K3$ surfaces},
J.\ Reine Angew.\ Math.\  588  (2005), 99--148.



\bibitem{freitag}
E.~Freitag:
\textit{A graded algebra related to cubic surfaces},
Kyushu J.\ Math.\  56  (2002),  no. 2, 299--312.


\bibitem{freitag2}
E.~Freitag:
\textit{Comparison of different models of the moduli
space of marked cubic surfaces}, Proc.\ Japanese-German
seminar Riushi-do, T.~Ibukyama, W.~Kohnen eds., 74-79 (2002),
also available from
\texttt{www.rzuser.uni-heidelberg.de/$\sim$t91/index4.html}.

\bibitem{geemen}
B.~van Geemen:
\textit{A linear system on Naruki's moduli space of marked cubic surfaces,}
Internat. J. Math. 13 (2002), 183-208.

\bibitem{hunt}
B.~Hunt:
\textit{The geometry of some special arithmetic
quotients}. Lecture Notes in
Mathematics, 1637, Springer-Verlag, Berlin (1996).


\bibitem{naruki}
 I.~Naruki:
\textit{Cross ratio variety as a moduli space of cubic surfaces},
Proc.\ London Math.\ Soc.\ 45 (1982), 1-30.

\bibitem{macdonald}
I.G.~Macdonald:
\textit{Some irreducible representations of Weyl groups},
Bull. London Math.\ Soc.\ 4 (1972), 148--150.

\bibitem{manin}
Y.~Manin:
\textit{Cubic Forms}, $2^{\text{nd}}$ edition,
North Holland Math.\ Library 4, North Holland 1986.


\bibitem{yoshida}
M.~Yoshida:
\textit{A $W(E_6)$-equivariant projective embedding of the moduli
space of cubic surfaces},
15 p.\, available at \texttt{arXiv:math.AG/0002102}

\bibitem{wall:git}
C.T.C.~Wall:
\textit{Geometric invariant theory of linear systems},
Math.\ Proc.\ Cambridge Philos.\ Soc.\  93  (1983), 57--62.

\bibitem{wall:nets}
C.T.C.~Wall:
\textit{Nets of quadrics, and theta-characteristics of singular curves},
Philos.\ Trans.\ Roy.\ Soc.\ London Ser. A  289  (1978), no. 1357, 229--269.

\end{thebibliography}
\end{document}